\documentclass[12pt]{article}

\usepackage{amssymb,latexsym, theorem}
\usepackage{hyperref}

\usepackage{graphicx}
\usepackage{pgf,tikz}
\usetikzlibrary{arrows}
\usepackage{color}

\def\d{\mathop{\mbox{d}}\nolimits}
\def\pr{\noindent{\bf Proof. }}
\def\eop{{\tiny \hfill $\blacksquare$}}

\def\N{{\mathbb N}}

\def\F{{\mathbb F}}
\def\O{{\mathbb O}}
\def\L{{\mathbb L}}

\def\PG{\mbox{\rm PG}}

\newcommand{\m}{\mathcal}
\newcommand{\mr}{\mathrm}
\newcommand{\dist}{\mathrm{d}}

\newtheorem{lemma}{Lemma}[section]
\newtheorem{theo}[lemma]{Theorem}

\newtheorem{co}[lemma]{Corollary}
\newtheorem{prop}[lemma]{Proposition}

\begin{document}

\title{The $\mathrm{L}_3(4)$ near octagon}
\author{A. Bishnoi and B. De Bruyn}
\maketitle

\begin{abstract}
In recent work we constructed two new near octagons, one related to the finite simple group $\mathrm{G}_2(4)$ and another one as a sub-near-octagon of the former. In the present paper, we give a direct construction of this sub-near-octagon using a split extension of the group $\mathrm{L}_3(4)$. 
We derive several geometric properties of this $\mathrm{L}_3(4)$ near octagon, and determine its full automorphism group. We also prove that the $\mathrm{L}_3(4)$ near octagon is closely related to the second subconstituent of the distance-regular graph on 486 vertices discovered by Soicher in 1993.
\end{abstract}


\section{Introduction} \label{sec1}

The main goal of this paper is to give a new construction of the near octagon of order $(2,4)$ that was discovered by us in \cite{ab-bdb:1}, and to derive some geometric and group theoretic properties of this near octagon that were not explored earlier. In doing so, we introduce a family of near polygons to which this new near octagon belongs and study some common properties that the members of this family have. This also gives an alternate treatment of some of the results in \cite[Section 4.1]{ab-bdb:1}. We emphasize that we have singled out a few combinatorial properties of the collinearity graph (properties (P1), (P2), (P3), (P4) in Section \ref{sec3}) that allows us to derive all the remaining properties. 

A connected partial linear space $\mathcal{S}=(\mathcal{P},\mathcal{L},\mathrm{I})$ with nonempty point set $\mathcal{P}$, line set $\mathcal{L}$ and incidence relation $\mathrm{I} \subseteq \mathcal{P} \times \mathcal{L}$ is called a {\em near polygon} if for every point $x$ and every line $L$, there exists a unique point on $L$ nearest to $x$. Here, distances are measured in the collinearity graph $\Gamma$. If $d$ is the diameter of $\Gamma$, then the near polygon is called a {\em near $2d$-gon}. A near polygon is said to have {\em order} $(s,t)$ if every line is incident with precisely $s+1$ points and every point is incident with exactly $t+1$ lines. 
For other relevant definitions see Section \ref{sec:prelim}. 

In \cite{ab-bdb:1}, we constructed a new near octagon related to the simple group $\mathrm{G}_2(4)$. Among other things, we proved the following.

\begin{prop} \label{prop1.1}
\begin{itemize}
\item[$(a)$] Let $\O_1$ be the point-line geometry whose points are the $4095$ central involutions of the group $G_1 = \mathrm{G}_2(4){:}2$ \footnote{This group is isomorphic to the automorphism group of $\mathrm{G}_2(4)$.} and whose lines are all the three element subsets $\{x, y, xy\}$ where $x$ and $y$ are two commuting central involutions that satisfy $[G_1:N_{G_1}(\langle x,y \rangle)] \in \{1365, 13650\}$, with incidence being  containment. Then $\O_1$ is a near octagon of order $(2,10)$.

\item[$(b)$] Let $S_1$ denote the set of all lines $\{ x,y,xy \}$ where $x$ and $y$ are two commuting central involutions of $G_1$ that satisfy $[G_1:N_{G_1}(\langle x,y \rangle)] =1365$. Then $S_1$ is a line spread of $\O_1$. If $\mathcal{Q}_1$ denotes the set of all quads of $\O_1$, then for every $Q \in \mathcal{Q}_1$ the lines of $S_1$ contained in $Q$ define a line spread of $Q$. Moreover, the point-line geometry $\mathcal{S}_1$ with point set $S_1$, line set $\mathcal{Q}_1$ and the natural incidence is a generalized hexagon isomorphic to the dual split Cayley hexagon $\mathrm{H}(4)^D$.
\end{itemize}
\end{prop}

\medskip \noindent In \cite{ab-bdb:1}, we also constructed a sub-near-octagon of $\O_1$ corresponding to each subhexagon of order $(4,1)$ of $\mathcal{S}_1 \cong \mathrm{H}(4)^D$. 
These (isomorphic) near octagons had order $(2,4)$ and they turned out to be new as well. 
We now give an alternate direct construction of this near octagon of order $(2,4)$ using the finite simple group $\mathrm{PSL}_3(4) = \mathrm{L}_3(4)$. The group $\mathrm{Aut}(\mathrm{L}_3(4)) \cong \mathrm{L}_3(4){:}D_{12}$ has a unique conjugacy class of subgroups of type $\mathrm{L}_3(4){:}2^2$, and we denote by $\mathrm{L}_3(4){:}2^2$ any isomorphic copy of such a subgroup.

\begin{theo} \label{theo1.2}
\begin{itemize}
\item[$(a)$] Let $\O_2$ be the point-line geometry whose points are the $315$ central involutions of the group $G_2 = \mathrm{L}_3(4){:}2^2$ and whose lines are all the three element subsets $\{x, y, xy\}$ where $x$ and $y$ are two commuting central involutions that satisfy $[G_2:N_{G_2}(\langle x,y \rangle)]\in \{ 105,420 \}$, with incidence being containment. Then $\O_2$ is a near octagon of order $(2,4)$.
\item[$(b)$] Let $S_2$ denote the set of all lines $\{ x,y,xy \}$ where $x$ and $y$ are two commuting central involutions of $G_2$ that satisfy $[G_2:N_{G_2}(\langle x,y \rangle)] =105$. Then $S_2$ is a line spread of $\O_2$. If $\mathcal{Q}_2$ denotes the set of all quads of $\O_2$, then for every $Q \in \mathcal{Q}_2$ the lines of $S_2$ contained in $Q$ define a line spread of $Q$. Moreover, the point-line geometry $\mathcal{S}_2$ with point set $S_2$, line set $\mathcal{Q}_2$ and natural incidence  is isomorphic to the unique generalized hexagon $\mathrm{H}(4,1)$ of order $(4,1)$.
\end{itemize}
\end{theo}

\bigskip \noindent We will refer to the near octagons $\O_1$ and $\O_2$  as the {\em $\mathrm{G}_2(4)$ near octagon} and the {\em $\mathrm{L}_3(4)$ near octagon}, respectively. 
We note that the central involutions of $\mathrm{L}_3(4){:}2^2$ are precisely the central involutions of its derived subgroup $\mathrm{L}_3(4)$, and these are precisely the nontrivial elations of the projective plane $\PG(2,4)$ on which the group $\mathrm{L}_3(4)$ acts naturally.

\medskip \noindent Part $(a)$ of Theorem \ref{theo1.2} will be proved in Section \ref{sec2}. The structure of the near octagons $\O_1$ and $\O_2$ around a fixed point can be described by a diagram, see \cite[Figure 1]{ab-bdb:1} for the $\mathrm{G}_2(4)$ near octagon and Figure \ref{suborbit} of the present paper for the $\mathrm{L}_3(4)$ near octagon. 
These two diagrams are very similar and Section \ref{sec3}  is devoted to the study of the geometric properties of a family of near octagons whose local structures can be described by such diagrams. 
Part $(b)$ of Theorem \ref{theo1.2} will follow from that discussion. 
We also note that Part $(b)$ of Proposition \ref{prop1.1} (with  the exception of the isomorphism $\mathcal{S}_1 \cong \mathrm{H}(4)^D$) also follows from the results of Section \ref{sec3}. 
At this moment, we do not know whether the $\mathrm{G}_2(4)$ and $\mathrm{L}_3(4)$ near octagons are the only nontrivial members of the family of near octagons under consideration in Section \ref{sec3}.  
In Section \ref{sec4}, we determine the full automorphism group of the near octagon $\O_2$. 

\begin{theo} \label{theo1.3}
We have $\mathrm{Aut}(\O_2) \cong \mathrm{L}_3(4) : 2^2$, and the automorphisms of $\O_2$ are precisely the conjugations by elements of $\mathrm{L}_3(4){:}2^2$.
\end{theo}

\bigskip \noindent  In \cite{So}, Soicher constructed a distance-transitive graph $\Upsilon$ of diameter $4$, which is a distance-regular graph with intersection array $\{ b_0,b_1,b_2,b_3;$ $c_1,c_2,c_3,c_4 \} = \{  56,45,16,1;$ $1,8,45,56 \}$. Let $G$ be the automorphism group of this graph $\Upsilon$. By \cite{So}, we know that if $x$ is a vertex of $\Upsilon$, then the stabilizer $G_x$ of $G$ with respect to $x$ is isomorphic to $\mathrm{L}_3(4){:}2^2$. 
Let $\O_2(x)$ denote the near octagon defined on the central involutions of $G_x \cong \mathrm{L}_3(4){:}2^2$ (in the sense of Theorem \ref{theo1.2}).

If $x$ is a vertex of $\Upsilon$, then by \cite{So}, we also know that $|\Upsilon_2(x)|=315$ and that the local graph $\Upsilon_x$ is isomorphic to the well-known Gewirtz graph, which is the unique strongly regular graph with parameters $(v,k,\lambda,\mu) = (56,10,0,2)$ \cite{Br-Ha}. 
We show in Section \ref{sec5} that the second subconstituent of $\Upsilon$ is closely related to the $\mathrm{L}_3(4)$ near octagon. Among other things, we prove the following there.

\begin{theo} \label{theo1.4}
Let $x$ be a vertex of $\Upsilon$. Then the following hold:
\begin{enumerate}
\item[$(1)$] For every vertex $y \in \Upsilon_2(x)$, the element-wise stabilizer of $\Upsilon_1(x) \cap \Upsilon_1(y)$ inside $G_x$ has order $2$ and is generated by a central involution $\sigma_y$ of $G_x \cong \mathrm{L}_3(4){:}2^2$. Moreover, the map $\theta: y \mapsto \sigma_y$ defines a bijection between $\Upsilon_2(x)$ and the set of central involutions of $G_x$.

\item[$(2)$] Let $\Gamma$ denote the graph defined on the set $\Upsilon_2(x)$ of vertices at distance $2$ from $x$ in $\Upsilon$, by making two vertices $y_1, y_2$ adjacent if and only if $\dist(y_1, y_2) = 4$ or ($\dist(y_1, y_2) = 2$ and $|\Upsilon_1(x) \cap \Upsilon_1(y_1) \cap \Upsilon_1(y_2)| = 4$). Then the map $\theta$ is an isomorphism between $\Gamma$ and the collinearity graph of $\O_2(x)$.

\item[$(3)$] The map $\theta$ also defines an isomorphism between the subgraph of $\Upsilon$ induced on $\Upsilon_2(x)$ (the second subconstituent) and the graph obtained from the collinearity graph of $\O_2(x)$ by making two vertices (points of the near octagon) adjacent when they are at distance $2$ from each other and have a unique common neighbour. \footnote{This graph is also a distance-regular graph, and recently Soicher \cite{Soicher16} has proved that it is the unique distance-regular graph with intersection array $\{32, 27, 8, 1; 1, 4, 27, 32\}$.}
\end{enumerate}
\end{theo}

\bigskip \noindent It is known that the second subconstituent of the Hall-Janko graph (which is a strongly regular graph $srg(100,36,14,12)$) is isomorphic to the distance-2 graph of the dual split Cayley hexagon $\mathrm{H}(2)^D$, and that the second subconstituent of the $\mathrm{G}_2(4)$-graph (which is an $srg(416,100,36,20)$) is isomorphic to the distance-2 graph of the Hall-Janko near octagon. One of the referees of an earlier paper of ours \cite{ab-bdb:2} informed us that there is also a close connection between the $\mathrm{G}_2(4)$ near octagon $\O_1$ and the second subconstituent of the distance-regular graph with intersection array $\{ 416,315,64,1;1,32,315,416 \}$ discovered by Soicher \cite{So} (see \cite[Remark 1.4]{ab-bdb:2}) 
We elaborate on this in the appendix of the present paper where we prove a similar result as Theorem \ref{theo1.4} for the $\mathrm{G}_2(4)$ near octagon.

\section{Preliminaries} \label{sec:prelim}

We first review some of the basics from the theory of near polygons that will be used in this paper. For more details see the standard reference \cite{bdb_book}.  

\bigskip \noindent Recall that a point-line geometry $\m S = (\m P, \m L, \mr I)$ is called a near $2d$-gon if it is a partial linear space with  collinearity graph of diameter $d$ and the property that for every  $x \in \m P$ and every  $L \in \m L$ there exists a unique point on $L$ that is nearest to $x$. 
Let $\m S = (\m P, \m L, \mr I)$ be a near $2d$-gon. 
Then by \cite[Lemma 1]{Br-Wi}, only one of the following possibilities occur for two lines $L_1$ and $L_2$ in $\m S$:
\begin{itemize}
\item there exists a unique point $x_1 \in L_1$ and a unique point $x_2 \in L_2$ such that for all $y_1 \in L_1$ and $y_2 \in L_2$, we have $\dist(y_1, y_2) = \dist(y_1, x_1) + \dist(x_1, x_2) + \dist(x_2, y_2)$;
\item for every point $y_1 \in L_1$, there exists a unique point $y_2 \in L_2$ such that $\dist(y_1, y_2) = \dist(L_1, L_2)$ and conversely for every $z_2 \in L_2$ there exists a unique point $z_1 \in L_1$ such that $\dist(z_1, z_2) = \dist(L_1, L_2)$. 
\end{itemize}
If the latter case occurs, then we say that the lines $L_1$ and $L_2$ are \textit{parallel}. 
A subset $S$ of $\m L$ is called a \textit{line spread} if for every point $x$, there exists a unique $L \in S$ such that $x \in L$. 

A subset $X$ of $\m P$ is called a \textit{subspace} if for every two collinear points $x, y$ in $X$, all points on the line $xy$ of $\m S$ are contained in $X$. For a nonempty subspace $X$ of $\m S$, we will denote by $\widetilde{X}$ the point-line geometry formed by taking the points of $X$ and those lines of $\m S$ which are completely contained in $X$. A subspace $X$ is called \textit{convex} if for every pair of points $x, y$ in $X$ every point on every shortest path between $x$ and $y$ in $\m S$ is also contained in $X$. A convex subspace $Q$ of $\m P$ is called a \textit{quad} of $\m S$ if $\widetilde{Q}$ is a non-degenerate generalized quadrangle (see \cite{Pa-Th} for the definition and basic properties of generalized quadrangles). From \cite[Proposition 2.5]{Sh-Ya}, it follows that if $\m S$ has at least three points on each line (which will be the case for all near polygons considered in this paper), then for every two points $x, y$ in $\m S$ at distance $2$ from each other which have more than one common neighbour, there exists a unique quad containing $x$ and $y$. 
From this it also follows that through a pair of intersecting lines in $\m S$ there can be at most one quad. 
If every line of $\m S$ has $s + 1 \geq 3$ points on it, and a pair of points at distance $2$ have $\alpha + 1 \geq 2$ common neighbours, then the generalized quadrangle corresponding to the unique quad through these points has order $(s, \alpha)$. 
From \cite[Proposition 2.6]{Sh-Ya}, one of the following cases occurs for a point-quad pair $(x, Q)$ in $\m S$: 
\begin{itemize}
\item there exists a unique point $x' \in Q$ such that for every $y \in Q$ we have $\dist(x, y) = \dist(x, x') + \dist(x', y)$;
\item the points in $Q$ nearest to $x$ form an \textit{ovoid} of the corresponding generalized quadrangle $\widetilde{Q}$, i.e., a set of points of $\widetilde{Q}$ meeting every line of $\widetilde{Q}$ in a singleton. 
\end{itemize}
In the former case we say that the point $x$ is \textit{classical} with respect to $Q$ and the point $x'$ is denoted by $\pi_Q(x)$. 
In the latter case, we say that $x$ is \textit{ovoidal} with respect to $Q$. 
If $\dist(x, Q) \leq 1$, then it is easy to show that $x$ is always classical with respect to $Q$. 

\bigskip \noindent 
We follow the group theoretical notation of Atlas \cite{Atlas}. An involution $\sigma$ of a group $G$ is called \textit{central} if there exists a Sylow-$2$ subgroup $H$ of $G$ such that $\sigma$ lies in the centralizer $C_G(H)$. 
We describe here a model of the group $\mathrm{L}_3(4){:}2^2$ which will be useful. 

Let $V$ be a 3-dimensional vector space over the field $\F_4$ with basis $(\bar e_1,\bar e_2,\bar e_3)$, and denote by $(\bar f_1,\bar f_2,\bar f_3)$ the dual basis of the dual space $V'$. The annihilator of a 1-space of $V'$ is a 2-space of $V$, implying that we can identify each line of $\PG(V)$ with a point of $\PG(V')$.  

With each \textit{collineation} $\theta \in P\Gamma L(V)$ of the projective space $\PG(V)$, there is associated a nonsingular $3 \times 3$ matrix $A$ over $\F_4$ and an automorphism $\tau$ of $\F_4$ such that the point  $\langle x_1 \bar e_1 + x_2 \bar e_2 + x_3 \bar e_3 \rangle$ of $\PG(V)$ is mapped to the point $\langle x_1' \bar e_1 + x_2' \bar e_2 + x_3' \bar e_3 \rangle$ of $\PG(V)$, where $[ x_1' \  x_2' \  x_3' ]^T = A \cdot [ x_1^\tau \  x_2^\tau \  x_3^\tau ]^T$. While the automorphism $\tau$ of $\F_4$ is uniquely determined by $\theta$, the matrix $A$ itself is only determined up to a nonzero factor. However, all matrices $A$ corresponding to $\theta$ have the same determinant as $k^3=1$ for every $k \in \F_4 \setminus \{ 0 \}$. The element $\theta$ of $P\Gamma L(V)$ also permutes the lines of $\PG(V)$. Specifically, if the ``line'' $\langle y_1 \bar f_1 + y_2 \bar f_2 + y_3 \bar f_3 \rangle$ is mapped to the ``line'' $\langle y_1' \bar f_1 + y_2' \bar f_2 + y_3' \bar f_3 \rangle$, then $[y_1' \ y_2' \ y_3']^T = (A^T)^{-1} \cdot [y_1^\tau \ y_2^\tau \ y_3^\tau]^T$.

With each \textit{correlation} $\theta$ of $\PG(V)$, there is also associated a nonsingular $3 \times 3$ matrix $A$ over $\F_4$ and an automorphism $\tau$ of $\F_4$ such that the point  $\langle x_1 \bar e_1 + x_2 \bar e_2 + x_3 \bar e_3 \rangle$ is mapped to the ``line'' $\langle y_1' \bar f_1 + y_2' \bar f_2 + y_3' \bar f_3 \rangle$, where $[y_1' \ y_2' \ y_3']^T = A \cdot [x_1^\tau \ x_2^\tau \ x_3^\tau]^T$. Again, the automorphism $\tau$ of $\F_4$  is uniquely determined by $\theta$, but $A$ is only determined up to a nonzero factor. If $A$ and $\tau$ are as above, then $\theta$ will map the ``line'' $\langle y_1 \bar f_1 + y_2 \bar f_2 + y_3 \bar f_3 \rangle$ to the point $\langle x_1' \bar e_1 + x_2' \bar e_2 + x_3' \bar e_3 \rangle$, where $[x_1' \  x_2'  \  x_3']^T = (A^T)^{-1} \cdot [y_1^\tau \ y_2^\tau \ y_3^\tau]^T$.

Let $G_1 \cong PS\mathrm{L}_3(4) = \mathrm{L}_3(4)$ denote the group of all collineations of $\PG(V)$ for which $\tau=1$ and for which $\det(A) = 1$,\footnote{This is well-defined. There are different choices possible for $A$, but the determinants are always the same.} let $G_2$ denote the group of all collineations of $\PG(V)$ for which $\det(A)=1$ and let $G$ denote the group of all collineations and correlations of $\PG(V)$ for which $\det(A)=1$. 
Then $G_1$ has index 2 in $G_2$ which itself has index 2 in $G$. If $G_3$ denotes the set of all elements of $G$ for which $A=I_3$, then $G_3 \cong C_2 \times C_2$ and $G$ is the internal semidirect product $G_1 \rtimes G_3$. So, $G$ is a group of type $\mathrm{L}_3(4){:}2^2$. Note also that $G_1$ is the derived subgroup of $G$, and that the group of all collineations and correlations of $\PG(V)$ has type $\mathrm{L}_3(4){:}D_{12}$. 

The group $G$ has a natural action on the set $\mathcal{F}$ of 105 flags of $\PG(V)$, i.e., the edges of the incidence graph of the projective plane $\PG(V)$. This action is transitive and there are four orbitals, implying that the group $G$ acts distance-transitively on the point set of the unique generalized hexagon $\mathrm{H}(4,1)$ of order $(4,1)$. This is the generalized hexagon whose points are the elements of $\mathcal{F}$ and whose lines are the points and lines of $\PG(V)$, with incidence being reverse containment. The fact that $G$ acts distance-transitively implies that $G$ must also act primitively on the set of 105 points of $\mathrm{H}(4,1)$. Consulting GAP's library of primitive permutation groups \cite{Gap}, we see that there exists a unique primitive permutation group on 105 letters that has type $\mathrm{L}_3(4){:}2^2$. We have used the implementation of this group in GAP to verify several claims in this paper (see \cite{ab-bdb:3}).

\section{The $\mathrm{L}_3(4)$ near octagon} \label{sec2}

\medskip \noindent 
Let $G \cong \mathrm{L}_3(4){:}2^2$ be the group of all collineations and correlations of $\mathrm{PG}(2, 4)$ whose corresponding matrix has determinant 1, as described in Section \ref{sec:prelim}. 
Let $\mathcal{P}$ denote the set of all 315 central involutions of $G$ (Type 2A in Atlas notations \cite{Atlas}). These $315$ central involutions are all contained in $G_1 = G' \cong \mathrm{L}_3(4)$ and correspond to the $315$ nontrivial \textit{elations} of the projective plane $\PG(2,4)$. 
The group $G$ acts on $\mathcal{P}$ by conjugation. Let $\omega$ denote a fixed central involution of $G$. 
Then it is easily checked using GAP that the stabilizer $G_\omega$ of $\omega$ has eight orbits on $\mathcal{P}$, the so-called {\em suborbits of} $G$ {\em with respect to} $\omega$. Such a suborbit will be denoted by $\mathcal{O}_0$, $\mathcal{O}_{1a}$, $\mathcal{O}_{1b}$, $\mathcal{O}_{2a}$, $\mathcal{O}_{2b}$, $\mathcal{O}_{3a}$, $\mathcal{O}_{3b}$ and $\mathcal{O}_4$ in accordance with the information provided by the following table:

\begin{center}
\begin{tabular}{|c||c|c|c|c|c|c|c|c|}
\hline
Suborbit & $\mathcal{O}_0$ & $\mathcal{O}_{1a}$ & $\mathcal{O}_{1b}$ & $\mathcal{O}_{2a}$ & $\mathcal{O}_{2b}$ & $\mathcal{O}_{3a}$ & $\mathcal{O}_{3b}$ & $\mathcal{O}_4$ \\
\hline
\hline
Size & 1 & 2 & 8 & 16 & 32 & 64 & 64 & 128 \\
\hline
$\langle x,\omega \rangle$ & $C_2$ & $C_2 \times C_2$ & $C_2 \times C_2$ & $C_2 \times C_2$ & $D_8$ & $D_8$ & $S_3$ & $D_{10}$ \\
\hline
\end{tabular}
\end{center}

\bigskip \noindent In the table, we have mentioned the sizes of the suborbits and descriptions for the groups $\langle x,\omega \rangle$, where $x$ is an arbitrary element of the considered suborbit. 
The suborbits with respect to a fixed central involution $y$ will be denoted by $\mathcal{O}_0(y),\mathcal{O}_{1a}(y),\ldots,\mathcal{O}_4(y)$.

The central involutions distinct from $\omega$ which commute with $\omega$ are those of the set $\mathcal{O}_{1a} \cup \mathcal{O}_{1b} \cup \mathcal{O}_{2a}$. Moreover, if $\mathcal{O} \in \{ \mathcal{O}_{1a},\mathcal{O}_{1b},\mathcal{O}_{2a} \}$ and $x \in \mathcal{O}$, then $\omega x \in \mathcal{O}$. The sets $\mathcal{O}_{1a}$, $\mathcal{O}_{1b}$ and $\mathcal{O}_{2a}$ consist of those central involutions $x \not= \omega$ for which $[G:N_G(\langle x,\omega \rangle)]$ has size 105, 420 and 840, respectively.

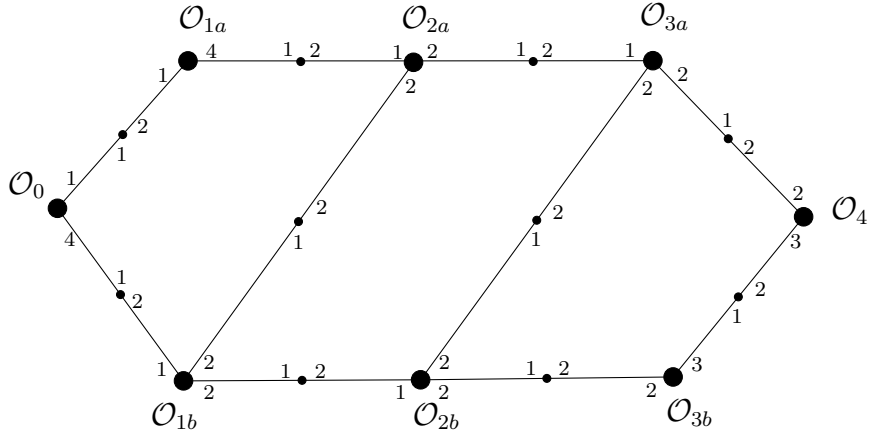
\begin{figure}
\begin{center}
\begin{tikzpicture}[line cap=round,line join=round,>=triangle 45,x=1.0cm,y=1.0cm]
\draw (0.3393083836366923,3.0573010907233384)-- (-1.39563107885606,1.0960651766010965);
\draw (-1.39563107885606,1.0960651766010965)-- (0.28,-1.2);
\draw (0.28,-1.2)-- (3.356594405363218,3.0573010907233384);
\draw (0.3393083836366923,3.0573010907233384)-- (3.356594405363218,3.0573010907233384);
\draw (3.356594405363218,3.0573010907233384)-- (6.52,3.06);
\draw (6.52,3.06)-- (8.523696717569893,0.9829169507863518);
\draw (0.28,-1.2)-- (3.432026555906381,-1.1857573773295889);
\draw (3.432026555906381,-1.1857573773295889)-- (6.52,3.06);
\draw (3.432026555906381,-1.1857573773295889)-- (6.788757255077141,-1.1480413020580074);
\draw (6.788757255077141,-1.1480413020580074)-- (8.523696717569893,0.9829169507863518);
\draw (3.03295311787782,3.933501017569001) node[anchor=north west] {$\mathcal{O}_{2a}$};
\draw (6.20,3.9669953715641197) node[anchor=north west] {$\mathcal{O}_{3a}$};
\draw (8.743740474045596,1.3879301139399627) node[anchor=north west] {$\mathcal{O}_4$};
\draw (6.55,-1.30) node[anchor=north west] {$\mathcal{O}_{3b}$};
\draw (3.183677710855855,-1.35) node[anchor=north west] {$\mathcal{O}_{2b}$};
\draw (-0.30,-1.35) node[anchor=north west] {$\mathcal{O}_{1b}$};
\draw (0.06870278930979522,3.933501017569001) node[anchor=north west] {$\mathcal{O}_{1a}$};
\draw (-2.2,1.722873653891152) node[anchor=north west] {$\mathcal{O}_0$};
\draw (1.45,3.45) node[anchor=north west] {{\scriptsize 1}};
\draw (4.55,3.45) node[anchor=north west] {{\scriptsize 1}};
\draw (7.29,2.50) node[anchor=north west] {{\scriptsize 1}};
\draw (7.42,-0.11) node[anchor=north west] {{\scriptsize 1}};
\draw (4.75,0.92) node[anchor=north west] {{\scriptsize 1}};
\draw (1.59,0.88) node[anchor=north west] {{\scriptsize 1}};
\draw (-0.76,2.05) node[anchor=north west] {{\scriptsize 1}};
\draw (-0.76,0.42) node[anchor=north west] {{\scriptsize 1}};
\draw (1.42,-0.83) node[anchor=north west] {{\scriptsize 1}};
\draw (4.68,-0.80) node[anchor=north west] {{\scriptsize 1}};
\draw (-0.48,2.42) node[anchor=north west] {{\scriptsize 2}}; % line
\draw (1.81,3.44) node[anchor=north west] {{\scriptsize 2}}; % line
\draw (4.90,3.43) node[anchor=north west] {{\scriptsize 2}}; % line
\draw (7.58,2.15) node[anchor=north west] {{\scriptsize 2}}; % line
\draw (-0.55,0.10) node[anchor=north west] {{\scriptsize 2}}; % line
\draw (1.90,1.35) node[anchor=north west] {{\scriptsize 2}}; % line
\draw (5.04,1.28) node[anchor=north west] {{\scriptsize 2}}; % line
\draw (1.88,-0.83) node[anchor=north west] {{\scriptsize 2}}; % line
\draw (5.12,-0.80) node[anchor=north west] {{\scriptsize 2}}; % line
\draw (7.73,0.25) node[anchor=north west] {{\scriptsize 2}}; % line
\draw (-1.45,0.92) node[anchor=north west] {{\scriptsize 4}}; % 0, t
\draw (-1.43,1.73) node[anchor=north west] {{\scriptsize 1}}; % 0, 1
\draw (-0.21,3.10) node[anchor=north west] {{\scriptsize 1}}; % 1', 1
\draw (0.43,3.42) node[anchor=north west] {{\scriptsize 4}}; % 1', t
\draw (2.91,3.41) node[anchor=north west] {{\scriptsize 1}}; % 2', 1
\draw (3.37,3.42) node[anchor=north west] {{\scriptsize 2}}; % 2', t - t'
\draw (3.08,2.95) node[anchor=north west] {{\scriptsize 2}}; % 2', t'
\draw (6.23,2.93) node[anchor=north west] {{\scriptsize 2}}; % 3', t'
\draw (8.20,0.89) node[anchor=north west] {{\scriptsize 3}}; % 4, t + t - t/t'
\draw (8.23,1.53) node[anchor=north west] {{\scriptsize 2}}; % 4, t/t'
\draw (6.70,3.11) node[anchor=north west] {{\scriptsize 2}}; % 3', t - t'
\draw (6.00,3.43) node[anchor=north west] {{\scriptsize 1}}; % 3', 1
\draw (6.89,-0.74) node[anchor=north west] {{\scriptsize 3}}; % 3'', t + 1 - t/t'
\draw (6.28,-1.09) node[anchor=north west] {{\scriptsize 2}}; % 3'', t/t'
\draw (3.52,-1.11) node[anchor=north west] {{\scriptsize 2}}; % 2'', t - t'
\draw (3.53,-0.72) node[anchor=north west] {{\scriptsize 2}}; % 2'', t'
\draw (2.95,-1.11) node[anchor=north west] {{\scriptsize 1}}; % 2'', 1
\draw (-0.23,-0.81) node[anchor=north west] {{\scriptsize 1}}; % 1'', 1
\draw (0.40,-0.71) node[anchor=north west] {{\scriptsize 2}}; % 1'', t'
\draw (0.40,-1.12) node[anchor=north west] {{\scriptsize 2}}; % 1'', t - t'
\begin{scriptsize}
\draw [fill=black] (-1.39563107885606,1.0960651766010965) circle (3.5pt);
\draw [fill=black] (0.3393083836366923,3.0573010907233384) circle (3.5pt);
\draw [fill=black] (0.28,-1.2) circle (3.5pt);
\draw [fill=black] (3.356594405363218,3.0573010907233384) circle (1.5pt);
\draw [fill=black] (6.52,3.06) circle (3.5pt);
\draw [fill=black] (8.523696717569893,0.9829169507863518) circle (3.5pt);
\draw [fill=black] (3.432026555906381,-1.1857573773295889) circle (3.5pt);
\draw [fill=black] (6.788757255077141,-1.1480413020580074) circle (3.5pt);
\draw [fill=black] (3.337736367727427,3.0384430530875477) circle (3.5pt);
\draw [fill=black] (1.8385223756820597,3.0478720719054433) circle (1.5pt);
\draw [fill=black] (4.928868183863713,3.049221526543774) circle (1.5pt);
\draw [fill=black] (7.521848358784946,2.021458475393176) circle (1.5pt);
\draw [fill=black] (7.656226986323517,-0.0825621756358278) circle (1.5pt);
\draw [fill=black] (5.1103919054917615,-1.166899339693798) circle (1.5pt);
\draw [fill=black] (4.97601327795319,0.9371213113352056) circle (1.5pt);
\draw [fill=black] (1.8560132779531906,-1.1928786886647944) circle (1.5pt);
\draw [fill=black] (1.8088681838637135,0.9192215265437739) circle (1.5pt);
\draw [fill=black] (-0.55781553942803,-0.05196741169945174) circle (1.5pt);
\draw [fill=black] (-0.5281613476096838,2.0766831336622174) circle (1.5pt);
\end{scriptsize}
\end{tikzpicture}
\end{center}
\caption{The suborbit diagram for the central involutions of $\mathrm{L}_3(4){:}2^2$} \label{suborbit}
\end{figure}

Let $\mathcal{L}$ denote the set of all triples $\{ x,y,xy \}$, where $x$ and $y$ are two distinct commuting central involutions such that $[G:N_G(\langle x,y \rangle)]$ has size 105 or 420, and let $\O_2$ be the point-line geometry with point set $\mathcal{P}$ and line set $\mathcal{L}$, where incidence is containment. With GAP (\cite{ab-bdb:3}, similar code as in \cite{ab-bdb:1}), we have computed the suborbit diagram for the central involutions of $\mathrm{L}_3(4){:}2^2$. This suborbit diagram can be found in Figure \ref{suborbit}. Each of the eight big nodes of the diagram denotes a suborbit and an edge between two such nodes denotes that there is a line of $\O_2$ that intersects both suborbits. A smaller node on each edge denotes a line and the two numbers on the node denote the number of points of the line that lie in each suborbit it intersects. The numbers on each big node are the number of lines through a representative of the suborbit meeting the adjacent suborbits.

\begin{theo} \label{theo2.1}
$\O_2$ is a near octagon of order $(2,4)$.
\end{theo}
\pr
Let $x$ be an arbitrary point of $\O_2$. For $i \in \N$ denote by $\Gamma_i(x)$,  the set of points at distance $i$ from $x$. Considering the suborbit diagram with respect to $x$, we see that:
\begin{itemize}
\item There are five lines through $x$.
\item $\Gamma_1(x)=\mathcal{O}_{1a}(x) \cup \mathcal{O}_{1b}(x)$, $\Gamma_2(x)=\mathcal{O}_{2a}(x) \cup \mathcal{O}_{2b}(x)$, $\Gamma_3(x)=\mathcal{O}_{3a}(x) \cup \mathcal{O}_{3b}(x)$ and $\Gamma_4(x) = \mathcal{O}_4(x)$. 
\item Every line contains a unique point nearest to $x$. 
\end{itemize}
So, $\O_2$ must be a near octagon of order $(2,4)$.
\eop

\bigskip \noindent For every point $x$ of $\O_2$, let $L_x$ denote the unique line through $x$ containing the two points of $\mathcal{O}_{1a}(x)$. Let $S_2$ denote the set of all lines of the form $L_x$, where $x \in \mathcal{P}$, and let $\mathcal{Q}$ denote the set of all quads of $\O_2$. 

\begin{lemma} \label{lem2.2}
For every point $x$ of $\O_2$, there are precisely two quads through $x$ and these two quads intersect in the line $L_x$.
\end{lemma}
\pr
We will rely on the information provided by the suborbit diagram.

If $y \in \mathcal{O}_{2b}(x)$, then $x$ and $y$ have precisely one neighbour and so they cannot be contained in a common quad. If $y \in \mathcal{O}_{2a}(x)$, then $x$ and $y$ have precisely three common neighbours one of which is contained on the line $L_x$, implying that $x$ and $y$ are contained in a unique quad of order $(2,2)$ (necessarily isomorphic to $W(2)$ \cite[5.2.3]{Pa-Th}) and this quad contains the line $L_x$.  

Every point $z \in \mathcal{O}_{1b}(x)$ is collinear with a point of $\mathcal{O}_{2a}(x)$, implying that the line $xz$ is contained in a quad. This quad must be unique as it also contains the line $L_x$. As there are four lines through $x$ distinct from $L_x$, there are precisely $\frac{4}{2}=2$ quads through $x$, and these two quads necessarily intersect in the line $L_x$.
\eop

\bigskip \noindent The following is an immediate consequence of Lemma \ref{lem2.2}.

\begin{co} \label{co2.3}
\begin{enumerate}
\item[$(1)$] The lines of $S_2$ are precisely the lines of $\O_2$ that are contained in two quads.
\item[$(2)$] The set $S_2$ is a line spread of $\O_2$.
\end{enumerate}
\end{co}

\bigskip \noindent In the following theorem, we show that there exists a full embedding of $\O_2$ into $\O_1$. In \cite{ab-bdb:1}, we constructed near octagons of order $(2,4)$ as subgeometries of $\O_1$. In Section \ref{sec3}, we explain why all these near octagons of order $(2,4)$ are isomorphic to the $\mathrm{L}_3(4)$ near octagon $\O_2$.

\begin{theo} \label{theo2.4}
There exists a full embedding of $\O_2$ into $\O_1$ mapping lines of $S_2$ to lines of $S_1$. Two embedded points are collinear in $\O_1$ if and only if they are collinear in $\O_2$. 
\end{theo}
\pr
By the Atlas \cite{Atlas} (also see \cite[Section 4.3.6]{Wi}), the group $G := \mathrm{G}_2(4){:}2$ has a maximal subgroup $H$ of type $\mathrm{SL}(3,4){:}2^2$. If we regard $G$ as naturally acting on the dual split Cayley hexagon $\mathrm{H}(4)^D$, then $H$ is the stabilizer of a subhexagon of order $(4,1)$. Let $Z := Z(H')$ be the center of the derived subgroup $H' \cong \mathrm{SL}(3,4)$ of $H$. Then $Z \cong C_3$, $J := H/Z \cong \mathrm{L}_3(4){:}2^2$ and $J' = H'/Z \cong \mathrm{L}_3(4)$. If $x \in H \cong \mathrm{SL}(3,4):2^2$, then $\overline{x}$ denotes the element $x Z$ of $H/Z \cong \mathrm{L}_3(4){:}2^2$.

There exists a bijective correspondence between the involutions of the group $H' \cong \mathrm{SL}(3,4)$ and the involutions of the group $J' = H'/Z \cong \mathrm{L}_3(4)$. If $\sigma$ is an involution of $H'$, then $\sigma \not\in Z \cong C_3$ and so $\overline{\sigma} = \sigma Z$ is an involution of $J' = H'/Z \cong \mathrm{L}_3(4)$. Conversely, if $\tau Z$ is an involution of $J' = H'/Z \cong \mathrm{L}_3(4)$, there exists a unique involution $\sigma$ of $H'$ for which $\overline{\sigma} = \tau Z$. Now, denote by $\Sigma$ the set of all involutions of $H'$ for which $\overline{\Sigma} := \{   \overline{\sigma} \, | \, \sigma \in \Sigma \}$ is the set of all central involutions of $J' = H'/Z \cong \mathrm{L}_3(4)$. Then $\Sigma$ is a conjugacy class of involutions of both $H'$ and $H$. If $\sigma \in \Sigma$ and $x \in H$, then $\overline{x^{-1} \sigma x} = \overline{x}^{-1} \cdot\overline{\sigma} \cdot \overline{x} \in \overline{\Sigma}$ as $\overline{x} \in H/Z \cong \mathrm{L}_3(4){:}2^2$ and $\overline{\sigma} \in \overline{\Sigma}$ and hence the involution $x^{-1} \sigma x$ belongs to $\Sigma$. 

As $\Sigma$ is a conjugacy class of involutions of $H'$ and $H$, it is contained in a conjugacy class of involutions of $G \cong \mathrm{G}_2(4){:}2$. We have verified with GAP (see \cite{ab-bdb:3}) that all elements of $\Sigma$ are central involutions of $G \cong \mathrm{G}_2(4){:}2$. Since $\overline{x^{-1} \sigma x} = \overline{x}^{-1} \cdot\overline{\sigma} \cdot \overline{x}$ for all $\sigma \in \Sigma$ and all $x \in H \cong \mathrm{SL}(3,4):2^2$, the action of $H \cong \mathrm{SL}(3,4):2^2$ on the set $\Sigma$ gives rise to the same suborbit diagram as the action of $H/Z \cong \mathrm{L}_3(4){:}2^2$ on $\overline{\Sigma}$ (which was depicted in Figure \ref{suborbit}). We suppose that $\O_2$ is the near octagon of order $(2,4)$ defined on the involutions of the set $\Sigma$, where two distinct involutions $\sigma_1$ and $\sigma_2$ of $\Sigma$ are collinear whenever the suborbit with respect to $\sigma_1$ that contains $\sigma_2$ has size 2 or 8.  We suppose that $\O_1$ is the near octagon of order $(2,10)$ defined on the central involutions of $G = \mathrm{G}_2(4){:}2$. Since all elements of $\Sigma$ are central involutions of $\mathrm{G}_2(4){:}2$, we have identified each point of $\O_2$ with a point of $\O_1$.

Now, fix a certain $\omega \in \Sigma$. We denote by $\mathcal{O}_0$, $\mathcal{O}_{1a}$, $\mathcal{O}_{1b}$, $\mathcal{O}_{2a}$, $\mathcal{O}_{2b}$, $\mathcal{O}_{3a}$, $\mathcal{O}_{3b}$, $\mathcal{O}_4$ the suborbits with respect to $\omega$ for the action of $H$ on $\Sigma$ (using the same notational convention as before). By \cite{ab-bdb:1}, we know that the action of $G \cong \mathrm{G}_2(4){:}2$ on the set of central involutions gives rise to a similar suborbit diagram as in Figure \ref{suborbit} (but with different numbers). We denote by $\mathcal{O}_0'$, $\mathcal{O}_{1a}'$, $\mathcal{O}_{1b}'$, $\mathcal{O}_{2a}'$, $\mathcal{O}_{2b}'$, $\mathcal{O}_{3a}'$, $\mathcal{O}_{3b}'$, $\mathcal{O}_4'$ the suborbits with respect to $\omega$ (using the same notational convention as in \cite{ab-bdb:1}). With the aid of GAP (see \cite{ab-bdb:3}), we have verified that every $\mathcal{O} \in \{ \mathcal{O}_0,\mathcal{O}_{1a},\mathcal{O}_{1b},\mathcal{O}_{2a},\mathcal{O}_{2b},\mathcal{O}_{3a},\mathcal{O}_{3b},\mathcal{O}_4 \}$ is contained in $\mathcal{O}'$.

We know that the suborbits $\mathcal{O}_{1a}$ and $\mathcal{O}_{1b}$ contain those points of $\O_2$ that are at distance 1 from $\omega$ (in $\O_2$). By \cite{ab-bdb:1}, we know that the suborbits $\mathcal{O}_{1a}'$ and $\mathcal{O}_{1b}'$ contain those points of $\O_1$ that are at distance 1 from $\omega$ (in $\O_1$). So, if $x$ and $y$ are two points of $\O_2$, then $x$ and $y$ are collinear in $\O_2$ if and only if $x$ and $y$ are collinear in $\O_1$. Every line of $\O_2$ is thus a collection of three mutually collinear points of $\O_1$, i.e. is a line of $\O_1$. 

The unique line of $S_2$ through $\omega$ is $\{ \omega \} \cup \mathcal{O}_{1a}$, while the unique line of $S_1$ containing $\omega$ is  $\{ \omega \} \cup \mathcal{O}_{1a}'$. As $\mathcal{O}_{1a} = \mathcal{O}_{1a}'$, every line of $S_2$ is also a line of $S_1$.
\eop

\section{A family of near octagons} \label{sec3}

Suppose $\mathcal{S}=(\mathcal{P},\mathcal{L},\mathrm{I})$ is a finite near octagon of order $(s,t)$, $s \geq 2$, and $S$ is a line spread of $\mathcal{S}$. For every $x \in \m P$, let $L_x$ denote the unique line of the line spread $S$ containing $x$. For every two points $x$ and $y$ of $\mathcal{S}$, $\d(x,y)$ denotes the distance between $x$ and $y$ in the collinearity graph $\Gamma$ of $\mathcal{S}$ and $\Gamma_i(x)$ with $i \in \{ 0,1,2,3,4 \}$ denotes the set of points at distance $i$ from $x$. We define a number of additional sets of points of $\mathcal{S}$.
\begin{itemize}
\item For every point $x$ of $\mathcal{S}$, we define $\Gamma_1'(x) := L_x \setminus \{ x \}$ and $\Gamma_1''(x) := \Gamma_1(x) \setminus \Gamma_1'(x)$.
\item For every point $x$ of $\mathcal{S}$ and every $i \in \{ 2,3 \}$, $\Gamma_i'(x)$ denotes the set of points of $\Gamma_i(x)$ that are collinear with a point of $\Gamma_{i-1}'(x)$, and we put $\Gamma_i''(x) := \Gamma_i(x) \setminus \Gamma_i'(x)$.
\end{itemize}
Throughout this section, we assume that there exists a positive divisor $t' \not= t$ of $t$ such that the following hold for every point $x$ of $\mathcal{S}$:
\begin{enumerate}
\item[(P1)] Every point of $\Gamma_2'(x)$ is incident with $t'$ lines meeting $\Gamma_1''(x)$.
\item[(P2)] Every point of $\Gamma_2''(x)$ is incident with a unique line meeting $\Gamma_1''(x)$.
\item[(P3)] Every point of $\Gamma_3'(x)$ is incident with $t'$ lines meeting $\Gamma_2''(x)$.
\item[(P4)] Every point of $\Gamma_3''(x)$ is incident with $\frac{t}{t'}$ lines meeting $\Gamma_2''(x)$.
\end{enumerate}

\bigskip \noindent From the suborbit diagram \cite[Figure 1]{ab-bdb:1}, it follows that the $\mathrm{G}_2(4)$ near octagon $\O_1$ is an example of such a near octagon with $s = 2$, $t = 10$ and $t' = 2$, and from Figure \ref{suborbit} it follows that the $\mathrm{L}_3(4)$ near octagon $\O_2$ is another example with $s = 2$, $t = 4$ and $t' = 2$. When $t' = 1$, we will show at the end of this section that $\m S$ must be the direct product of a generalized hexagon of order $(s, t - 1)$ and a line of size $s + 1$. We call this a \textit{trivial example}. We do not know whether the above family of near octagons contains nontrivial examples besides $\O_1$ and $\O_2$.

In the lemmas below, $x$ denotes some fixed point of $\mathcal{S}$. A point $y$ of $\mathcal{S}$ is said to be of {\em type} $i \in \{ 0,1,2,3,4 \}$ if $y \in \Gamma_i(x)$. A point $y$ of $\mathcal{S}$ is said to be of {\em type} $i'$ or $i''$, for $i \in \{ 1,2,3 \}$, if it belongs to $\Gamma_i'(x)$ or $\Gamma_i''(x)$, respectively. A line $L$ of $\mathcal{S}$ is said to have {\em type} $(i_1,i_2)$ with $i_1,i_2 \in \{ 0,1,1',1'',2,2',2'',3,3',3'',4 \}$ if $L$ contains a unique point of type $i_1$ and $s$ points of type $i_2$.

\medskip \noindent The following lemma follows from easy counting. 

\begin{lemma} \label{lem3.1}
We have $|\Gamma_0(x)|=1$, $|\Gamma_1'(x)|=s$ and $|\Gamma_1''(x)|=st$. There is a unique line of type $(0,1')$, namely $L_x$, and every line of the line spread $S$ meeting $\Gamma_0(x) \cup \Gamma_1'(x)$ coincides with $L_x$. The point $x$ is contained in a unique line of type $(0,1')$ and $t$ lines of type $(0,1'')$. 
\end{lemma}

\begin{lemma} \label{lem3.2}
Every point $y \in \Gamma_2'(x)$ is contained in a (necessarily unique) quad together with $x$. This quad has order $(s,t')$.
\end{lemma}
\pr
By definition of the set $\Gamma_2'(x)$, the line $L_x = \{ x \} \cup \Gamma_1'(x)$ contains a (necessarily unique) point collinear with $y$. By Property (P1), there are exactly $t'$ points in $\Gamma_1''(x)$ collinear with $y$. As the points $x$ and $y$ have  $t'+1 \geq 2$ common neighbours, they are  contained in a unique quad, necessarily of order $(s,t')$.
\eop

\begin{lemma} \label{lem3.3}
There does not exist any $y \in \Gamma_2''(x)$ which is contained in a quad together with $x$.
\end{lemma}
\pr
By Property (P2), there is a unique line through $y$ meeting $\Gamma_1''(x)$, and by definition of the set $\Gamma_2''(x)$, there are no lines through $y$ meeting $\Gamma_1'(x)$. Hence, $x$ and $y$ have a unique common neighbour, implying that they cannot be contained in a quad.
\eop

\begin{lemma} \label{lem3.4}
Every quad through a point $y$ of $\mathcal{S}$ contains the line $L_y$.
\end{lemma}
\pr
Without loss of generality, we may take $y=x$. Suppose $Q$ is a quad through $x$ not containing the line $L_x$. Then no $y \in \Gamma_2(x) \cap Q$ is collinear with a point of $L_x$ as otherwise $L_x$ would be contained in $Q$. This implies that $\Gamma_2(x) \cap Q \subseteq \Gamma_2''(x)$, which contradicts Lemma \ref{lem3.3}.
\eop

\begin{lemma} \label{lem3.5}
For every quad $Q$ of $\mathcal{S}$, the lines of $S$ contained in $Q$ define a line spread of the corresponding generalized quadrangle $\widetilde{Q}$. 
\end{lemma}
\pr
If this were not the case, then there would exist a line of $S$ meeting $Q$ in a single point, which would contradict Lemma \ref{lem3.4}.
\eop

\begin{lemma} \label{lem3.6}
Through every point $y \in \Gamma_1'(x)$, there is a unique line of type $(0,1')$ and $t$ lines of type $(1',2')$. 
\end{lemma}
\pr
By the definition of the set $\Gamma_2'(x)$, every line through $y$ distinct from $L_x = xy$ must have type $(1',2')$. 
\eop

\begin{lemma} \label{lem3.7}
We have $|\Gamma_2'(x)|=s^2t$.
\end{lemma}
\pr
By Lemmas \ref{lem3.1} and \ref{lem3.6}, the number of edges between $\Gamma_1'(x)$ and $\Gamma_2'(x)$ is equal to $|\Gamma_1'(x)| \cdot t \cdot s = s^2t$. As any  point of $\Gamma_2'(x)$ is collinear with a unique point of $L_x$, we see that $|\Gamma_2'(x)|=s^2t$.
\eop

\begin{lemma} \label{lem3.8}
There are precisely $\frac{t}{t'}$ quads through any given point of $\mathcal{S}$.
\end{lemma}
\pr
Without loss of generality, we prove this for the point $x$. By Lemmas \ref{lem3.2} and \ref{lem3.3}, every quad through $x$ has order $(s,t')$ and contains $s^2t'$ points at distance 2 from $x$ which all belong to $\Gamma_2'(x)$. From Lemmas \ref{lem3.2} and \ref{lem3.7}, it then follows that there are precisely $\frac{|\Gamma_2'(x)|}{s^2t'}=\frac{s^2t}{s^2t'} = \frac{t}{t'}$ quads through $x$.
\eop

\begin{lemma} \label{lem3.9}
Every line $M$ through a point $y$ of $\mathcal{S}$ distinct from $L_y$ is contained in a unique quad. This quad also contains the line $L_y$.
\end{lemma}
\pr
Without loss of generality, we may suppose that $y=x$. As any quad through $x$ contains the line $L_x$ (Lemma \ref{lem3.4}), there is at most one quad through $M$. Since every quad through $x$ has order $(s,t')$, the $\frac{t}{t'}$ quads through $x$ cover $t' \cdot \frac{t}{t'} = t$ lines through $x$ distinct from $L_x$. As these are all lines through $x$ distinct from $L_x$, $M$ must be contained in a unique quad.
\eop

\begin{lemma} \label{lem3.10}
For every $y \in \Gamma_1''(x)$, the line $L_y$ has type $(1'',2')$. For every $z \in \Gamma_2'(x)$, the line $L_z$ has type $(1'',2')$.
\end{lemma}
\pr
By Lemma \ref{lem3.9}, the line $xy$ is contained in a unique quad $Q$. By Lemma \ref{lem3.5}, $L_x$ and $L_y$ are two disjoint lines contained in $Q$, and so every point of $L_y \setminus \{ y \}$ belongs to $\Gamma_2'(x)$, as it is collinear with a point of $L_x \setminus \{ x \} = \Gamma_1'(x)$. Therefore, $L_y$ has type $(1'',2')$. 

The $st$ mutually distinct lines $L_y$ with $y \in \Gamma_1''(x)$ cover $s^2t$ points of $\Gamma_2'(x)$. By Lemma \ref{lem3.7}, these are all the points of $\Gamma_2'(x)$. So, for every $z \in \Gamma_2'(x)$, we also know that the line $L_z$ has type $(1'',2')$.
\eop

\begin{lemma} \label{lem3.11}
Let $y \in \Gamma_1''(x)$. Then $y$ is contained in a unique line of type $(0,1'')$, $t'$ lines of type $(1'',2')$ and $t-t'$ lines of type $(1'',2'')$.
\end{lemma}
\pr
There is a unique line through $y$ containing $x$, namely the line $xy$, and this line contains precisely $s$ points of $\Gamma_1''(x)$. Every line through $y$ distinct from $xy$ contains a unique point collinear with $x$ (namely $y$) and hence precisely $s$ points of $\Gamma_2(x)$. 

By Lemma \ref{lem3.9}, there is a unique quad $Q$ (of order $(s,t')$) through the line $xy$. By Lemma \ref{lem3.3}, the $t'$ lines of $Q$ through $y$ distinct from $xy$ all contain precisely $s$ points of $\Gamma_2'(x)$. Conversely, suppose that $L$ is a line through $y$ containing a point $u \in \Gamma_2'(x)$. The unique quad through $x$ and $u$ is a convex subspace and contains therefore the lines $L=yu$ and $xy$, implying that this quad coincides with $Q$. So, $L$ is one of the $t'$ lines of $Q$ through $y$ distinct from $xy$. The remaining $t-t'$ lines through $y$ must be of type $(1'',2'')$.
\eop

\begin{lemma} \label{lem3.12}
We have $|\Gamma_2''(x)|=s^2t(t-t')$.
\end{lemma}
\pr
By Lemmas \ref{lem3.1} and \ref{lem3.11}, the number of edges between $\Gamma_1''(x)$ and $\Gamma_2''(x)$ is equal to $|\Gamma_1''(x)| \cdot (t-t') s = s^2t(t-t')$. By Property (P2), we know that the number of edges is also equal to $|\Gamma_2''(x)|$.
\eop

\begin{lemma} \label{lem3.13}
Let $y \in \Gamma_2'(x)$. Then $y$ is incident with a unique line of type $(1',2')$, $t'$ lines of type $(1'',2')$ and $t-t'$ lines of type $(2',3')$.
\end{lemma}
\pr
The lines through $y$ meeting $\Gamma_1(x)$ (necessarily in a unique point) are precisely the $t'+1$ lines through $y$ that are contained in the unique quad $Q$ (of order $(s,t')$) through $x$ and $y$. By Lemma \ref{lem3.3}, each of these lines contains precisely $s$ points of $\Gamma_2'(x)$ and hence they have type $(1'',2')$ or $(1',2')$. Note also that there is a unique line through $y$ meeting $L_x = \{  x \} \cup \Gamma_1'(x)$, which is the unique line of type $(1',2')$ through $y$.

By definition of the set $\Gamma_3'(x)$, the remaining $t-t'$ lines through $y$ all contain $s$ points of $\Gamma_3'(x)$ and thus have type $(2',3')$.
\eop

\begin{lemma} \label{lem3.14}
Through every point $y \in \Gamma_3'(x)$, there is a unique line meeting $\Gamma_2'(x)$.
\end{lemma}
\pr
By definition of the set $\Gamma_3'(x)$, we know that there is at least one such a line. Suppose two lines through $y$ meet $\Gamma_2'(x)$ in the respective points $u_1$ and $u_2$. The point $u_i$ with $i \in \{ 1,2 \}$ is collinear with a unique point $v_i \in L_x$. This point necessarily coincides with the unique point $v$ of $L_x$ at distance 2 from $y$, i.e., $v=v_1=v_2$. Since the points $y$ and $v$ have two distinct neighbours, namely $u_1$ and $u_2$, they are contained in a unique quad $Q$. The line $L_x$ meets this quad $Q$ in a unique point (namely $v$), but this contradicts Lemma \ref{lem3.5}.
\eop

\begin{lemma} \label{lem3.15}
We have $|\Gamma_3'(x)| = s^3t(t-t')$.
\end{lemma}
\pr
By Lemmas \ref{lem3.7} and \ref{lem3.13}, the number of edges between $\Gamma_2'(x)$ and $\Gamma_3'(x)$ is equal to $|\Gamma_2'(x)| \cdot (t-t')s = s^3t(t-t')$. By Lemma \ref{lem3.14}, the number of these edges is also equal to $|\Gamma_3'(x)|$. 
\eop

\begin{lemma} \label{lem3.16}
For every $y \in \Gamma_2''(x)$, the line $L_y$ has type $(2'',3')$. For every $z \in \Gamma_3'(x)$, the line $L_z$ has type $(2'',3')$.
\end{lemma}
\pr
By Property (P2), there exists a unique point $u \in \Gamma_1''(x)$ collinear with $y$. By Lemma \ref{lem3.10}, the line $L_u$ has type $(1'',2')$ and so is distinct from $uy$. By Lemma \ref{lem3.9}, there is a unique quad $Q$ through $L_u$ and $uy$. By Lemma \ref{lem3.5}, this quad also contains the line $L_y$. By Lemma \ref{lem3.3}, $Q$ cannot contain $x$. So, $x$ lies at distance 1 from $Q$ and is classical with respect to $Q$, implying that $L_y \setminus \{ y \} \subseteq \Gamma_3(x)$. The lines $L_y$ and $L_u$ of $Q$ are parallel and at distance 1. Hence, every point of $L_y \setminus \{ y \}$ is collinear with a point of $L_u \setminus \{ u \} \subseteq \Gamma_2'(x)$, implying that $L_y \setminus \{ y \} \subseteq \Gamma_3'(x)$. So, $L_y$ has type $(2'',3')$.

The $s^2t(t-t')$ lines $L_y$ with $y \in \Gamma_2''(x)$ all have type $(2'',3')$ and are mutually disjoint, implying that they cover $s^3t(t-t')$ points of 
$\Gamma_3'(x)$. By Lemma \ref{lem3.15}, these are all the points of $\Gamma_3'(x)$, and so the second claim of the lemma should be true as well.  
\eop

\begin{lemma} \label{lem3.18}
For any $y \in \Gamma_3''(x) \cup \Gamma_4(x)$, the line $L_y$ has type $(3'',4)$.
\end{lemma}
\pr
This follows from Lemmas \ref{lem3.1}, \ref{lem3.10}, \ref{lem3.16} and the fact that $S$ is a line spread of $\mathcal{S}$.
\eop

\begin{lemma} \label{lem3.19}
There are no lines meeting $\Gamma_3'(x)$ and $\Gamma_3''(x)$.
\end{lemma}
\pr
Suppose $L$ is a line containing points $u' \in \Gamma_3'(x)$ and $u'' \in \Gamma_3''(x)$, and let $y$ denote the unique point of $\Gamma_2(x)$ on $L$. By the definition of $\Gamma_3''(x)$, we must have $y \in \Gamma_2''(x)$. By Lemma \ref{lem3.14}, the point $u'$ is collinear with a unique point $v \in \Gamma_2'(x)$. By Lemma \ref{lem3.10}, the line $L_v$ meets $\Gamma_1''(x)$ and so is distinct from $vu'$. By Lemma \ref{lem3.9}, there exists a unique quad $Q$ through $L_v$ and $vu'$. 
This quad contains a point of $\Gamma_1''(x)$, and so the $t'$ lines of $Q$ through $u' \in \Gamma_3'(x)$ distinct from $vu'$ all meet $\Gamma_2(x)$. By Lemma \ref{lem3.14}, we know that these $t'$ lines meet $\Gamma_2''(x)$. By Property (P3), we know that these are all the lines through $u'$ meeting $\Gamma_2''(x)$. 
So, the line $L$ is one of these lines and hence contained in $Q$. But then every point of $L \setminus \{ y \}$ is collinear with a point of $L_v \setminus \Gamma_1''(x) \subseteq \Gamma_2'(x)$, implying that each point of $L \setminus \{ y \}$ belongs to $\Gamma_3'(x)$. This contradicts the fact that $u'' \in \Gamma_3''(x)$.  
\eop

\bigskip \noindent The following is a consequence of Lemmas \ref{lem3.1}, \ref{lem3.6}, \ref{lem3.10}, \ref{lem3.11}, \ref{lem3.13}, \ref{lem3.16}, \ref{lem3.18} and \ref{lem3.19}.

\begin{co} \label{co3.20}
\begin{itemize}
\item Every line of $\mathcal{S}$ has type $(0,1')$, $(0,1'')$, $(1',2')$, $(1'',2')$, $(1'',2'')$, $(2',3')$, $(2'',3')$, $(2'',3'')$, $(3',4)$ or $(3'',4)$.
\item Every line of the line spread $S$ has type $(0,1')$, $(1'',2')$, $(2'',3')$ or $(3'',4)$.
\end{itemize}
\end{co}

\begin{figure}
\begin{center}
\begin{tikzpicture}[line cap=round,line join=round,>=triangle 45,x=1.0cm,y=1.0cm]
\draw (0.3393083836366923,3.0573010907233384)-- (-1.39563107885606,1.0960651766010965);
\draw (-1.39563107885606,1.0960651766010965)-- (0.28,-1.2);
\draw (0.28,-1.2)-- (3.356594405363218,3.0573010907233384);
\draw (0.3393083836366923,3.0573010907233384)-- (3.356594405363218,3.0573010907233384);
\draw (3.356594405363218,3.0573010907233384)-- (6.52,3.06);
\draw (6.52,3.06)-- (8.523696717569893,0.9829169507863518);
\draw (0.28,-1.2)-- (3.432026555906381,-1.1857573773295889);
\draw (3.432026555906381,-1.1857573773295889)-- (6.52,3.06);
\draw (3.432026555906381,-1.1857573773295889)-- (6.788757255077141,-1.1480413020580074);
\draw (6.788757255077141,-1.1480413020580074)-- (8.523696717569893,0.9829169507863518);
\draw (3.03295311787782,3.933501017569001) node[anchor=north west] {$\Gamma_{2}'(x)$};
\draw (6.20,3.9669953715641197) node[anchor=north west] {$\Gamma_{3}'(x)$};
\draw (8.743740474045596,1.3879301139399627) node[anchor=north west] {$\Gamma_4(x)$};
\draw (6.55,-1.30) node[anchor=north west] {$\Gamma_{3}''(x)$};
\draw (3.183677710855855,-1.35) node[anchor=north west] {$\Gamma_{2}''(x)$};
\draw (-0.30,-1.35) node[anchor=north west] {$\Gamma_{1}''(x)$};
\draw (0.06870278930979522,3.933501017569001) node[anchor=north west] {$\Gamma_{1}'(x)$};
\draw (-2.85,1.412873653891152) node[anchor=north west] {$\Gamma_0(x)$};
\draw (1.45,3.45) node[anchor=north west] {{\scriptsize 1}};
\draw (4.55,3.45) node[anchor=north west] {{\scriptsize 1}};
\draw (7.29,2.50) node[anchor=north west] {{\scriptsize 1}};
\draw (7.42,-0.11) node[anchor=north west] {{\scriptsize 1}};
\draw (4.75,0.92) node[anchor=north west] {{\scriptsize 1}};
\draw (1.59,0.88) node[anchor=north west] {{\scriptsize 1}};
\draw (-0.76,2.05) node[anchor=north west] {{\scriptsize 1}};
\draw (-0.76,0.42) node[anchor=north west] {{\scriptsize 1}};
\draw (1.42,-0.83) node[anchor=north west] {{\scriptsize 1}};
\draw (4.68,-0.80) node[anchor=north west] {{\scriptsize 1}};
\draw (-0.48,2.42) node[anchor=north west] {{\scriptsize $s$}}; % line
\draw (1.81,3.44) node[anchor=north west] {{\scriptsize $s$}}; % line
\draw (4.90,3.43) node[anchor=north west] {{\scriptsize $s$}}; % line
\draw (7.58,2.15) node[anchor=north west] {{\scriptsize $s$}}; % line
\draw (-0.55,0.10) node[anchor=north west] {{\scriptsize $s$}}; % line
\draw (1.90,1.35) node[anchor=north west] {{\scriptsize $s$}}; % line
\draw (5.04,1.28) node[anchor=north west] {{\scriptsize $s$}}; % line
\draw (1.88,-0.83) node[anchor=north west] {{\scriptsize $s$}}; % line
\draw (5.12,-0.80) node[anchor=north west] {{\scriptsize $s$}}; % line
\draw (7.73,0.25) node[anchor=north west] {{\scriptsize $s$}}; % line
\draw (-1.45,0.92) node[anchor=north west] {{\scriptsize $t$}}; % 0, t
\draw (-1.43,1.73) node[anchor=north west] {{\scriptsize 1}}; % 0, 1
\draw (-0.21,3.10) node[anchor=north west] {{\scriptsize 1}}; % 1', 1
\draw (0.43,3.42) node[anchor=north west] {{\scriptsize $t$}}; % 1', t
\draw (2.91,3.41) node[anchor=north west] {{\scriptsize 1}}; % 2', 1
\draw (3.37,3.45) node[anchor=north west] {{\scriptsize $t - t'$}}; % 2', t - t'
\draw (3.08,2.95) node[anchor=north west] {{\scriptsize $t'$}}; % 2', t'
\draw (6.23,2.93) node[anchor=north west] {{\scriptsize $t'$}}; % 3', t'
\draw (8.20,0.89) node[anchor=north west] {{\scriptsize $t + 1 - \frac{t}{t'}$}}; % 4, t + t - t/t'
\draw (8.15,1.80) node[anchor=north west] {{\scriptsize $\frac{t}{t'}$}}; % 4, t/t'
\draw (6.70,3.11) node[anchor=north west] {{\scriptsize $t - t'$}}; % 3', t - t'
\draw (6.00,3.43) node[anchor=north west] {{\scriptsize 1}}; % 3', 1
\draw (6.89,-0.74) node[anchor=north west] {{\scriptsize $t + 1 - \frac{t}{t'}$}}; % 3'', t + 1 - t/t'
\draw (6.28,-1.09) node[anchor=north west] {{\scriptsize $\frac{t}{t'}$}}; % 3'', t/t'
\draw (3.52,-1.05) node[anchor=north west] {{\scriptsize $t - t'$}}; % 2'', t - t'
\draw (3.53,-0.72) node[anchor=north west] {{\scriptsize $t'$}}; % 2'', t'
\draw (2.95,-1.11) node[anchor=north west] {{\scriptsize 1}}; % 2'', 1
\draw (-0.23,-0.81) node[anchor=north west] {{\scriptsize 1}}; % 1'', 1
\draw (0.40,-0.71) node[anchor=north west] {{\scriptsize $t'$}}; % 1'', t'
\draw (0.40,-1.12) node[anchor=north west] {{\scriptsize $t - t'$}}; % 1'', t - t'
\begin{scriptsize}
\draw [fill=black] (-1.39563107885606,1.0960651766010965) circle (3.5pt);
\draw [fill=black] (0.3393083836366923,3.0573010907233384) circle (3.5pt);
\draw [fill=black] (0.28,-1.2) circle (3.5pt);
\draw [fill=black] (3.356594405363218,3.0573010907233384) circle (1.5pt);
\draw [fill=black] (6.52,3.06) circle (3.5pt);
\draw [fill=black] (8.523696717569893,0.9829169507863518) circle (3.5pt);
\draw [fill=black] (3.432026555906381,-1.1857573773295889) circle (3.5pt);
\draw [fill=black] (6.788757255077141,-1.1480413020580074) circle (3.5pt);
\draw [fill=black] (3.337736367727427,3.0384430530875477) circle (3.5pt);
\draw [fill=black] (1.8385223756820597,3.0478720719054433) circle (1.5pt);
\draw [fill=black] (4.928868183863713,3.049221526543774) circle (1.5pt);
\draw [fill=black] (7.521848358784946,2.021458475393176) circle (1.5pt);
\draw [fill=black] (7.656226986323517,-0.0825621756358278) circle (1.5pt);
\draw [fill=black] (5.1103919054917615,-1.166899339693798) circle (1.5pt);
\draw [fill=black] (4.97601327795319,0.9371213113352056) circle (1.5pt);
\draw [fill=black] (1.8560132779531906,-1.1928786886647944) circle (1.5pt);
\draw [fill=black] (1.8088681838637135,0.9192215265437739) circle (1.5pt);
\draw [fill=black] (-0.55781553942803,-0.05196741169945174) circle (1.5pt);
\draw [fill=black] (-0.5281613476096838,2.0766831336622174) circle (1.5pt);
\end{scriptsize}
\end{tikzpicture}
\end{center}
\caption{The structure of $\mathcal{S}$ with respect to a fixed point $x$}  \label{fig:sub1}
\end{figure}
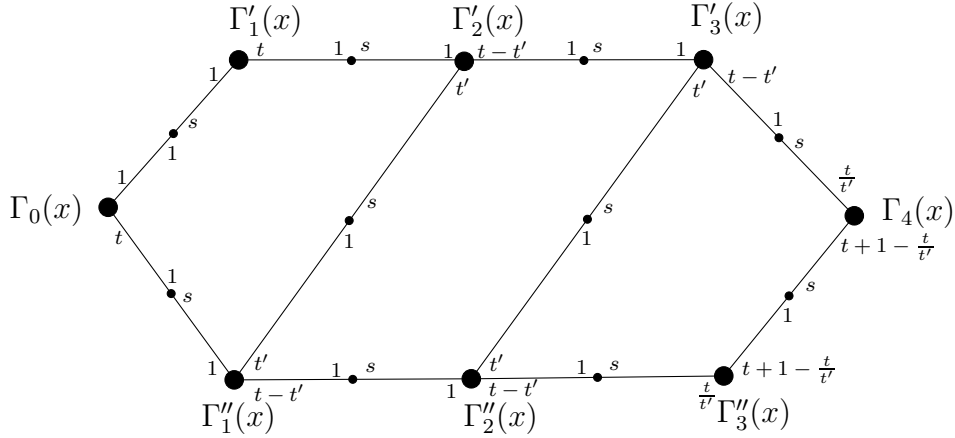

\bigskip \noindent We can now draw a diagram for $\mathcal{S}$ which looks similar to Figure \ref{suborbit}. This is done in Figure \ref{fig:sub1}, though we still need to prove some numerical information that is contained in the diagram.

\begin{lemma} \label{lem3.21}
Let $y \in \Gamma_2''(x)$, $z$ the unique point of $\Gamma_1''(x)$ collinear with $y$ and $Q$ the unique quad through the line $yz \not\in S$. Then the $t'$ lines of $Q$ through $y$ distinct from $yz$ have type $(2'',3')$.
\end{lemma}
\pr
Recall that by Property (P2), there is indeed a unique point $z \in \Gamma_1''(x)$ collinear with $y$. By Lemma \ref{lem3.10}, the line $L_z$ contains $s$ points of $\Gamma_2'(x)$. By Lemma \ref{lem3.5}, $L_z$ is contained in $Q$. If $L$ is one of the $t'$ lines of $Q$ through $y$ distinct from $yz$, then every point of $L \setminus \{ y \}$ is collinear with a point of $L_z \setminus \{ z \} \subseteq \Gamma_2'(x)$ and hence is contained in $\Gamma_3'(x)$ by the definition of $\Gamma_3'(x)$.
\eop

\begin{lemma} \label{lem3.22}
Every point $y$ of $\Gamma_3'(x)$ is incident with a unique line of type $(2',3')$, $t'$ lines of type $(2'',3')$ and $t-t'$ lines of type $(3',4)$.
\end{lemma}
\pr
By Lemma \ref{lem3.14}, there is a unique line of type $(2',3')$ through $y$ and by Property (P3), there are precisely $t'$ lines of type $(2'',3')$ through $y$. The remaining $t-t'$ lines through $y$ should all have type $(3',4)$.
\eop

\begin{lemma} \label{lem3.23}
Every point $y$ of $\Gamma_2''(x)$ is incident with a unique line of type $(1'',2'')$, $t'$ lines of type $(2'',3')$ and $t-t'$ lines of type $(2'',3'')$.
\end{lemma}
\pr
By Lemmas \ref{lem3.12} and \ref{lem3.21}, the number of edges between $\Gamma_2''(x)$ and $\Gamma_3'(x)$ is at least $|\Gamma_2''(x)| \cdot st' = s^3t t' (t-t')$, with equality if and only if every point of type $2''$ is incident with precisely $t'$ lines of type $(2'',3')$. By Lemma \ref{lem3.15} and Property (P3), we know that the number of edges between $\Gamma_2''(x)$ and $\Gamma_3'(x)$ is precisely $|\Gamma_3'(x)| \cdot t' = s^3 tt'(t-t')$. So, the point $y \in \Gamma_2''(x)$ is incident with precisely $t'$ lines of type $(2'',3')$. By Property (P2), $y$ is incident with a unique line of type $(1'',2'')$. The remaining $t-t'$ lines through $y$ should have type $(2'',3'')$.
\eop

\begin{lemma} \label{lem3.26}
Every point of $\Gamma_3''(x)$ is contained in $\frac{t}{t'}$ lines of type $(2'',3'')$ and $t+1-\frac{t}{t'}$ lines of type $(3'',4)$.
\end{lemma}
\pr
This follows from Property (P4)\footnote{Note that this is the first time we have used this property.} and the discussion so far regarding the types of lines in $\mathcal{S}$, where we deduced that the only types of lines through a point of type $3''$ are $(2'',3'')$ and $(3'',4)$.
\eop

\begin{lemma} \label{lem3.27}
We have $|\Gamma_3''(x)| = s^3t'(t-t')^2$.
\end{lemma}
\pr
By Lemmas \ref{lem3.12} and \ref{lem3.23}, the number of edges between $\Gamma_2''(x)$ and $\Gamma_3''(x)$ is equal to $|\Gamma_2''(x)| \cdot (t-t') s = s^3t(t-t')^2$. By Lemma \ref{lem3.26}, this number is also equal to $|\Gamma_3''(x)| \cdot \frac{t}{t'}$. Hence, $|\Gamma_3''(x)| = s^3t'(t-t')^2$.
\eop

\begin{lemma} \label{lem3.24}
Let $L$ be a line of $S$ meeting $\Gamma_3''(x)$ and $\Gamma_4(x)$. If $Q$ is a quad through $L$, then $Q$ contains at most one point at distance 2 from $x$.
\end{lemma}
\pr
Since the quad $Q$ contains a point at distance $4$ from $x$, it contains a point at distance $2$ from $x$ if and only if $x$ is classical with respect to $Q$. But then $\d(x,Q)=2$ and $Q$ contains a unique point at distance 2 from $x$.
\eop

\begin{lemma} \label{lem3.25}
Let $L$ be a line of $S$ having type $(3'',4)$. Then all $\frac{t}{t'}$ quads through $L$ contain a unique point of type $2''$ and a unique line of type $(2'',3')$. This point of type $2''$ is the unique point of the quad at distance $2$ from $x$, and this line of type $(2'',3')$ belongs to $S$.
\end{lemma}
\pr
Let $y$ be the unique point of $L$ contained in $\Gamma_3''(x)$. By Property (P4), there are $\frac{t}{t'}$ lines $K$ through $y$ meeting $\Gamma_2''(x)$ giving rise by Lemma \ref{lem3.24} to $\frac{t}{t'}$ distinct quads $\langle K,L \rangle$ through $L$. These are all the quads through $L$ by Lemma \ref{lem3.8}. Let $Q$ be any of these $\frac{t}{t'}$ quads and let $z$ denote the unique point in $Q \cap \Gamma_2(x)$. By Lemmas \ref{lem3.21} and \ref{lem3.23}, the $t'$ lines of type $(2'',3')$ through $z \in \Gamma_2''(x)$ are all contained in a quad $Q'$ through $z$. As $Q'$ contains a point of $\Gamma_1''(x)$, we have $Q \not= Q'$ and so $Q \cap Q'$ is at most a line. As $L_z \subseteq Q \cap Q'$, we have $L_z = Q \cap Q'$. Obviously, $L_z = Q \cap Q'$ is the unique line of type $(2'',3')$ of $Q$.
\eop

\begin{lemma} \label{lem3.28}
Every point $y$ of $\Gamma_4(x)$ is incident with $\frac{t}{t'}$ lines of type $(3',4)$ and $t+1-\frac{t}{t'}$ lines of type $(3'',4)$.
\end{lemma}
\pr 
If $L$ is a line through $y$ meeting $\Gamma_3'(x)$, then Lemmas \ref{lem3.9} and \ref{lem3.18} imply that there is a unique quad through $L$ and this quad contains the line $L_y$ of type $(3'',4)$. Conversely, by Lemma \ref{lem3.25}, any of the $\frac{t}{t'}$ quads through $L_y$ contains a unique line through $y$ meeting $\Gamma_3'(x)$. So, there are $\frac{t}{t'}$ lines through $y$ meeting $\Gamma_3'(x)$ and $t+1-\frac{t}{t'}$ lines meeting $\Gamma_3''(x)$.
\eop

\begin{lemma} \label{lem3.29}
We have $|\Gamma_4(x)| = s^4 t' (t-t')^2$.
\end{lemma}
\pr
By Lemmas \ref{lem3.15} and \ref{lem3.22}, the number of edges between $\Gamma_3'(x)$ and $\Gamma_4(x)$ is equal to $|\Gamma_3'(x)| \cdot (t-t')s = s^4t(t-t')^2$. By Lemma \ref{lem3.28}, the number of such edges is also equal to $|\Gamma_4(x)| \cdot \frac{t}{t'}$. It follows that $|\Gamma_4(x)| = s^4 t' (t-t')^2$.
\eop

\bigskip \noindent 
Thus we have determined all the numerical information given in Figure \ref{fig:sub1}. 
We now prove that the geometry defined on elements of the line spread $S$ and the set of quads of $\m S$ form a generalized hexagon. 

\begin{lemma} \label{lem3.30}
Every point-quad pair $(y,Q)$ is classical.
\end{lemma}
\pr
Without loss of generality, we may suppose that $y=x$. We have already remarked above that $(x,Q)$ is classical if $\d(x,Q) \leq 1$. By Lemmas \ref{lem3.4} and \ref{lem3.25}, $(x,Q)$ is classical if $Q$ contains points at distance 4 from $x$. So, suppose $Q \subseteq \Gamma_2(x) \cup \Gamma_3(x)$. Let $y \in Q \cap \Gamma_2(x)$. Then $L_y \subseteq Q$ implies by Lemma \ref{lem3.10} that $y \not\in \Gamma_2'(x)$. So, $y \in \Gamma_2''(x)$. The line $L_y$ meets $\Gamma_3'(x)$. By (P2), $y$ is collinear with a unique point $z$ of $\Gamma_1''(x)$. Denote by $Q'$ the unique quad through $yz$ and $L_y$. By Lemmas \ref{lem3.21} and \ref{lem3.23}, the $t'$ lines of type $(2'',3')$ through $y$ are the $t'$ lines of $Q'$ through $y$ distinct from $yz$. Each of the $t-t'$ lines through $y$ meeting $\Gamma_3''(x)$ determines a quad together with $L_y$ and such a quad has $t'$ lines through $y$ meeting $\Gamma_3''(x)$. So, there are $\frac{t-t'}{t'} = \frac{t}{t'}-1$ quads through $L_y$ meeting $\Gamma_3''(x)$ and hence also $\Gamma_4(x)$ as the unique line of $S$ through a point of $\Gamma_3''(x)$ meets $\Gamma_4(x)$. Taking into account $Q'$, we have thus accounted for all $\frac{t}{t'}$ quads through $y$, and none of these quads can be contained in $\Gamma_2(x) \cup \Gamma_3(x)$.
\eop

\begin{lemma} \label{lem3.31}
Two distinct lines $L_1$ and $L_2$ of $S$ lie at distance $1$ from each other if and only if they are contained in the same quad. If $\d(L_1,L_2)=1$, then $L_1$ and $L_2$ are parallel.
\end{lemma}
\pr
If $L_1$ and $L_2$ are contained in a quad, then they are parallel and $\d(L_1,L_2)=1$. Conversely, suppose that $\d(L_1,L_2)=1$ and let $x_1 \in L_1$ and $x_2 \in L_2$ such that $\d(x_1,x_2)=1$. The unique quad through the line $x_1x_2$ then contains $L_1$ and $L_2$.
\eop

\begin{lemma} \label{lem3.32}
Every two lines $K$ and $L$ of $S$ are parallel.
\end{lemma}
\pr
Let $k \in K$ and $l \in L$ such that $\d(K,L) = \d(k,l)$. We must show that every point of $K$ has distance at most $\delta := \d(k,l)$ from $L$. Let $k=y_0,y_1,\ldots,y_\delta=l$ be a path connecting $k$ and $l$. Put $L_i := L_{y_i}$, $i \in \{ 0,1,\ldots,\delta \}$. Then for every $i \in \{ 1,2,\ldots,\delta \}$, the lines $L_{i-1}$ and $L_i$ are either equal or parallel and at distance $1$. So, every point of $L_{i-1}$ has distance at most 1 from a point of $L_i$. This implies that every point of $K=L_0$ has distance at most $\delta$ from $L = L_\delta$.
\eop

\bigskip \noindent Let $\mathcal{Q}$ denote the set of quads of $\mathcal{S}$.
For the rest of this section let $\mathcal{S}'$ be the point-line geometry with point set $S$ and line set $\mathcal{Q}$, where incidence is containment.

\begin{lemma} \label{lem3.33}
If $L_1,L_2 \in S$, then the distance $d$ between $L_1$ and $L_2$ in $\mathcal{S}'$ is equal to $\d(L_1,L_2)$.
\end{lemma}
\pr
Put $d_1:=d$ and $d_2:=\d(L_1,L_2)$. Suppose $L_1=K_0,K_1,\ldots,K_{d_1}=L_2$ is a (shortest) path in $\mathcal{S}'$ connecting $L_1$ and $L_2$. Then the $K_i$'s are mutually disjoint. Let $x_0$ be an arbitrary point of $K_0$ and for every $i \in \{ 1,2,\ldots,d_1 \}$, let $x_i$ denote the unique point of $K_i$ collinear with $x_{i-1}$. Then $x_0 \in L_1$ lies at distance at most $d_1$ from $x_{d_1} \in L_2$, showing that $d_2 \leq d_1$.

Suppose $y_0,y_1,\ldots,y_{d_2}$ is a (shortest) path of length $d_2$ in $\mathcal{S}$ connecting a point $y_0 \in L_1$ with a point $y_{d_2} \in L_2$. In the sequence $L_{y_0},L_{y_1},\ldots,L_{y_{d_2}}$, any two consecutive lines are either equal or at distance 1 from each other (in $\mathcal{S}'$) by Lemma \ref{lem3.31}. So, we should also have that $d_1 \leq d_2$.

We conclude that $d_1=d_2$.
\eop

\begin{lemma} \label{lem3.34}
$\mathcal{S}'$ has order $(st',\frac{t}{t'}-1)$.
\end{lemma}
\pr
By Lemma \ref{lem3.5}, every line of $\mathcal{S}'$ contains precisely $1+st'$ points (namely the size of a line spread of a generalized quadrangle of order $(s,t')$). By Lemmas \ref{lem3.4} and \ref{lem3.8}, every line of $S$ is contained in precisely $\frac{t}{t'}$ quads, showing that every point of $\mathcal{S}'$ is incident with precisely $\frac{t}{t'}$ lines of $\mathcal{S}'$.
\eop

\begin{lemma} \label{lem3.35}
$\mathcal{S}'$ is a near hexagon.
\end{lemma}
\pr
Let $L \in S$ and $Q \in \mathcal{Q}$. Denote by $S_Q$ the set of lines of $S$ contained in $Q$. By Lemma \ref{lem3.5}, $S_Q$ is a line spread of the generalized quadrangle $\widetilde{Q}$. Let $x$ be an arbitrary point of $L$. As $x$ is classical with respect to the quad $Q$, there exists a unique point $x' \in Q$ nearest to $x$. By Lemma \ref{lem3.33}, $L_{x'}$ is the unique line of $S_Q$ at distance $\d(x,x')$ from $L_x=L$ in the geometry $\mathcal{S}'$. Every other line of $S_Q$ has distance $\d(x,x')+1$ from $L_x=L$.

It remains to show that $\mathcal{S}'$ has diameter 3. By Lemma \ref{lem3.33}, the diameter of $\mathcal{S}'$ is at most 3. If $x_1$ and $x_2$ are two points of $\mathcal{S}$ at distance 4 from each other, then the fact that $L_{x_1}$ and $L_{x_2}$ are parallel implies that they lie at distance 3 from each other (both in $\mathcal{S}$ as $\mathcal{S}'$). So, the diameter of $\mathcal{S}'$ is indeed 3.
\eop

\begin{lemma} \label{lem3.36}
Every two points of $\mathcal{S}'$ at distance $2$ from each other have a unique common neighbour.
\end{lemma}
\pr
Let $L_1,L_2 \in S$ be at distance 2 from each other in the geometry $\mathcal{S}'$ and suppose $M_1,M_2 \in S$ are two distinct neighbours of $L_1,L_2$. The lines $L_1,L_2,M_1,M_2$ are mutually parallel. By Lemma \ref{lem3.33}, there exist points $x_1 \in L_1$ and $x_2 \in L_2$ at distance $2$ from each other. By Lemmas \ref{lem3.32} and \ref{lem3.33}, the points $x_1,x_2$ have at least two common neighbours, one on $M_1$ and another one on $M_2$. So, $x_1$ and $x_2$ are contained in a quad $Q$. This quad should contain the line $L_{x_1}=L_1$, but that is impossible since $\d(x_2,L_1)=2$. 
\eop

\bigskip \noindent The following is an immediate consequence of Lemmas \ref{lem3.34}, \ref{lem3.35} and \ref{lem3.36}.

\begin{co} \label{co3.37}
$\mathcal{S}'$ is a generalized hexagon of order $(st',\frac{t}{t'}-1)$.
\end{co}

\bigskip \noindent Theorem \ref{theo1.2}(b) is an immediate consequence of Figure \ref{suborbit} and Corollary \ref{co3.37} (with $s=2$, $t'=2$ and $t=4$). Similarly, Proposition \ref{prop1.1}(b) (with exception of the isomorphism $\mathcal{S}_1 \cong \mathrm{H}(4)^D$) is a consequence of Corollary \ref{co3.37} and Figure 1 of \cite{ab-bdb:1} (with $s=2$, $t'=2$ and $t=10$). 

\medskip \noindent By Theorem \ref{theo2.4}, there exists a full embedding of $\O_2$ to $\O_1$ mapping $S_2$ of $S_1$. Via this full embedding, the geometry $\mathcal{S}_2 \cong \mathrm{H}(4,1)$ becomes a subgeometry of $\mathcal{S}_1 \cong \mathrm{H}(4)^D$, showing that the near octagons of order $(2,4)$ constructed in \cite[Theorem 1.2]{ab-bdb:1} are isomorphic to $\O_2$. 

\bigskip \noindent The finite near octagons $\mathcal{S}=(\mathcal{P},\mathcal{L},\mathrm{I})$ of order $(s,t)$ with $s \geq 2$ having a line spread $S$ for which Properties (P1), (P2), (P3) and (P4) are satisfied with $t'=1$ can be completely described. Suppose $\mathcal{S}$ is the direct or cartesian product $\mathcal{H} \times \L_{s+1}$ of a generalized hexagon $\mathcal{H}$ of order $(s,t-1)$ and a line $\mathbb{L}_{s+1}$ of size $s+1 \geq 3$ (as defined in Brouwer and Wilbrink \cite[Section (a)]{Br-Wi}). Then there exist two partitions $S$ and $\mathcal{U}$ of $\mathcal{P}$ such that the following hold:
\begin{itemize}
\item $S$ is a line spread.
\item Every $U \in \mathcal{U}$ is a convex subspace and $\widetilde{U} \cong \mathcal{H}$.
\item Every line of $S$ intersects every element of $\mathcal{U}$ in a singleton.
\end{itemize}
Then $(\mathcal{S},S)$ satisfies the Properties (P1) -- (P4) with $t'=1$. Moreover, if $x \in \mathcal{P}$ and $U_x$ denotes the unique element of $\mathcal{U}$ containing $x$, then
\begin{itemize}
\item $\Gamma_i''(x)$ with $i \in \{ 1,2,3 \}$ is the set of points of $U_x$ at distance $i$ from $x$;
\item $\Gamma_i'(x)$ with $i \in \{ 1,2,3 \}$ is the set of points of $\mathcal{P} \setminus U_x$ at distance $i$ from $x$. 
\end{itemize}
We refer to the above example as a {\em trivial member} of the family of near octagons. We can now prove the following.

\begin{prop} \label{prop3.38}
Suppose $(\mathcal{S},S)$ satisfies the conditions (P1) -- (P4) with $t'=1$. Then $\mathcal{S}$ is the cartesian product of a generalized hexagon of order $(s,t-1)$ and a line of size $s+1$. 
\end{prop}
\pr
Fix a line $L^\ast \in S$. For every $x \in L^\ast$, let $H_x$ denote the set of all points $y$ of $\mathcal{S}$ for which $x$ is the unique point on $L^\ast$ nearest to $y$. As lines of $S$ are mutually parallel (see Lemma \ref{lem3.32}), we know that the following holds.
\begin{enumerate}
\item[(a)] Every line of $S$ contains a unique point of each $H_x$, $x \in L^\ast$.
\end{enumerate}
Now, let $Q$ be a quad of $\mathcal{S}$. As $Q$ has order $(s,t')=(s,1)$, $\widetilde{Q}$ is an $(s+1) \times (s+1)$-grid. The set $S_Q$ of lines of $S$ contained in $Q$ is a line spread of $\widetilde{Q}$. Let $S_Q'$ denote the set of all remaining lines of $\widetilde{Q}$. Then $S_Q'$ is the second line spread of $\widetilde{Q}$. If $L$ is the unique line of $S_Q$ nearest to $L^\ast$ in $\mathcal{S}'$, then every point $x \in L^\ast$ is classical with respect to $Q$ with $\pi_Q(x) \in L$, implying the following:
\begin{enumerate}
\item[(b)] A quad $Q$ of $\mathcal{S}$ intersects each $H_x$, $x \in L^\ast$, in a line of $S_Q'$.
\end{enumerate}
We also know that every line of $\mathcal{S}$ not belonging to $S$ is contained in a unique quad.

By (a) and (b), we now see that each $H_x$, $x \in L^\ast$, is a subspace for which $\widetilde{H_x} \cong \mathcal{S}'$ is a generalized hexagon of order $(s,\frac{t}{t'}-1)=(s,t-1)$ and that $\mathcal{S} \cong \widetilde{H_x} \times \mathbb{L}_{s+1}$.
\eop

\section{Determination of $\mathrm{Aut}(\O_2)$} \label{sec4}

This section is devoted to the determination of the full automorphism group $\mathrm{Aut}(\O_2)$ of the near octagon $\O_2$. As before, we denote by $\mathcal{Q}_2$ the set of all quads (isomorphic to $W(2)$) of $\O_2$. The set of all lines of $\O_2$ that are contained in two quads is then a line spread $S_2$ of $\O_2$. Recall that the point-line geometry $\mathcal{S}_2$ defined by $S_2$ and $\mathcal{Q}_2$ is a generalized hexagon of order $(4,1)$. For every point $x$ of $\O_2$, we denote by $L_x$ the unique line of $S_2$ containing $x$.

\medskip \noindent The following lemma is well known and straightforward to prove.

\begin{lemma} \label{lem4.1}
\begin{enumerate}
\item[$(1)$] If $\sigma$ is a nontrivial elation with center $x$ and axis $L$ and $\theta$ is a collineation of $\PG(2,4)$, then $\theta^{-1} \sigma \theta$ is a nontrivial elation with center $x^\theta$ and axis $L^\theta$.
\item[$(2)$] If $\sigma$ is a nontrivial elation with center $x$ and axis $L$ and $\theta$ is a correlation of $\PG(2,4)$, then $\theta^{-1} \sigma \theta$ is a nontrivial elation with center $L^\theta$ and axis $x^\theta$.
\end{enumerate}
\end{lemma}

\begin{lemma} \label{lem4.2}
Every automorphism $\theta$ of $\O_2$ permutes the elements of the line spread $S_2$ and hence determines an automorphism of the generalized hexagon $\mathcal{S}_2 \cong \mathrm{H}(4,1)$ defined by $S_2$ and $\mathcal{Q}_2$.
\end{lemma}
\pr
The automorphism $\theta$ permutes the quads of $\O_2$ and hence the lines of $\O_2$ that can be obtained as intersection of two quads (these are the elements of $S_2$).
\eop

\begin{lemma} \label{lem4.3}
Let $Q$ be a quad of $\O_2$ and let $L_1,L_2,L_3,L_4,L_5$ the five lines of $S_2$ contained in $Q$. If $\theta$ is an automorphism of $\O_2$ fixing each of these lines, then $\theta$ fixes each point of $Q$.
\end{lemma}
\pr
Let $x$ be an arbitrary point of $Q$, $L = \{ x,y,z \}$ a line through $x$ not belonging to $S_2$ and contained in $Q$. 
Since $\widetilde{Q}$ is isomorphic to $W(2)$, the line $L$ is the unique line of $Q$ meeting $L_x$, $L_y$ and $L_z$. As $L_x^\theta = L_x$, $L_y^\theta = L_y$ and $L_z^\theta = L_z$, it follows that $L^\theta=L$ and $x^\theta = x$.
\eop

\bigskip \noindent The following is a consequence of Lemma \ref{lem4.3}.

\begin{co} \label{co4.5}
If $\theta$ is an automorphism of $\O_2$ fixing each line of $S_2$, then $\theta$ is trivial.
\end{co}

\begin{lemma} \label{lem4.4}
Let $Q$ be a quad of $\O_2$, and $L_1,L_2,L_3,L_4,L_5$ the five lines of $S_2$ contained in $Q$. If $\theta$ is an automorphism of $\O_2$ such that $L_1^\theta = L_1$, $L_2^\theta = L_2$, $L_3^\theta = L_4$, $L_4^\theta = L_5$ and $L_5^\theta = L_3$, then $\theta$ permutes the points of $L_1$ according to a cycle of length $3$.
\end{lemma}
\pr
For every $i \in \{ 3,4,5 \}$, choose the collinear points $a_i \in L_1$ and $b_i \in L_2$ such that $a_ib_i$ is the unique line of $Q$ meeting $L_1$, $L_2$ and $L_i$. As $L_1^\theta = L_1$, $L_2^\theta = L_2$, $L_3^\theta = L_4$, $L_4^\theta = L_5$ and $L_5^\theta = L_3$, we have $a_3^\theta = a_4$, $a_4^\theta = a_5$ and $a_5^\theta = a_3$.
\eop

\bigskip \noindent Recall from Lemma \ref{lem4.2} that every automorphism $\O_2$ induces an automorphism of $\mathcal{S}_2 \cong \mathrm{H}(4,1)$. By Corollary \ref{co4.5}, we then know the following.

\begin{co} \label{co4.6}
Every automorphism of $\mathcal{S}_2$ is induced by at most one automorphism of $\O_2$.
\end{co}

\bigskip \noindent Now, let $\widetilde{G}$ denote the subgroup of $\mathrm{Aut}(\mathcal{S}_2) \cong \mathrm{L}_3(4){:}D_{12}$ consisting of all automorphisms $\widetilde{\theta}$ of $\mathcal{S}_2$ that are induced by automorphisms of $\O_2$.

Recall that the points of $\mathcal{S}_2 \cong \mathrm{H}(4,1)$ can be identified with the flags of $\PG(2,4)$ and the lines of $\mathcal{S}_2$ with the points and lines of $\PG(2,4)$. Consider now a subgroup $H$ of $\mathrm{Aut}(\mathcal{S}_2)$ with order 3 corresponding to a group of homologies of $\PG(2,4)$ with axis $M$ and center $p \not\in M$.

\begin{lemma} \label{lem4.7}
We have $H \cap \widetilde{G} = \{ 1 \}$.
\end{lemma}
\pr
Suppose $\widetilde{\theta} \in H \cap \widetilde{G} \setminus \{ 1 \}$ is induced by an automorphism $\theta$ of $\O_2$. Let $Q_1$ be the quad of $\O_2$ corresponding to $M$, and let $Q_2$ be a quad of $\O_2$ corresponding to a point $q$ of $M$. Then $Q_1 \cap Q_2$ is the line of $S_2$ corresponding to the flag $(q,M)$ of $\PG(2,4)$. 

Now, every homology with center $p$ and axis $M$ fixes each point of $M$, implying that all lines of $S_2$ contained in $Q_1$ are fixed by $\theta$. By Lemma \ref{lem4.3}, this implies that $\theta$ fixes each point of $Q_1 \cap Q_2$.

Every nontrivial homology with center $p$ and axis $M$ fixes the lines $pq$, $M$ and permutes the other three lines through $q$ according to a cycle of length 3. This implies that $\theta$ fixes two lines of $S_2$ contained in $Q_2$ (including $Q_1 \cap Q_2$) and permutes the other three according to a cycle of length 3. By Lemma \ref{lem4.4}, this implies that $\theta$ permutes the points of $Q_1 \cap Q_2$ according to a cycle of length $3$. This is however impossible as we have already shown that $\theta$ fixes all points of $Q_1 \cap Q_2$.
\eop

\begin{prop} \label{prop4.8}
We have $\mathrm{Aut}(\O_2) \cong \mathrm{L}_3(4){:}2^2$. The automorphisms of $\O_2$ are given by conjugations by elements of $\mathrm{L}_3(4){:}2^2$.
\end{prop}
\pr
From Lemma \ref{lem4.1} it follows that the action of the group $\mathrm{L}_3(4){:}D_{12}$ and hence of the group $\mathrm{L}_3(4){:} 2^2$ on the elations via conjugation is faithful. Hence, $\mathrm{L}_3(4){:}2^2$ can be regarded as a group of automorphisms of $\O_2$, where the action is given by conjugation. To complete the proof of the proposition, it suffices to show that $|\mathrm{Aut}(\O_2)| \leq |\mathrm{L}_3(4){:}2^2|$. But that is not so hard if one invokes Lemma \ref{lem4.7}. As $\widetilde{G} \cap H = \{ 1 \}$, we have $3 \cdot |\widetilde{G}| = |\widetilde{G}| \cdot |H| = |\widetilde{G} H| \leq |\mathrm{L}_3(4){:}D_{12}|$, implying that $|\mathrm{Aut}(\O_2)| = |\widetilde{G}| \leq |\mathrm{L}_3(4){:}2^2|$.
\eop

\section{Connection with a distance-regular graph discovered by Soicher} \label{sec5}

Let $\Gamma$ denote the Gewirtz graph, that is the unique strongly regular graph with parameters $(v,k,\lambda,\mu)=(56,10,0,2)$. The automorphism group $\mathrm{Aut}(\Gamma)$ of $\Gamma$ is isomorphic to $\mathrm{L}_3(4){:}2^2$. We have computationally verified (see \cite{ab-bdb:3}) that if $\sigma$ is a central involution of $\mathrm{Aut}(\Gamma) \cong \mathrm{L}_3(4){:}2^2$, then the set $X_\sigma$ of vertices fixed by $\sigma$ is a set of size 8 and that the subgraph of $\Gamma$ induced on $X_\sigma$ is the union of two cycles of length 4. Moreover, the central involution $\sigma$ is uniquely determined by the set $X_\sigma$: the element-wise stabilizer of $X_\sigma$ inside $\mathrm{Aut}(\Gamma)$ is a group of order 2 generated by $\sigma$. Any set of vertices of $\Gamma$ of the form $X_\sigma$ for some central involution $\sigma$ of $\mathrm{Aut}(\Gamma)$ will be called a {\em special 8-set} of $\Gamma$. 

In \cite{So}, L. Soicher constructed a distance-regular  graph $\Upsilon$ (which is also distance-transitive) of diameter 4 with intersection array $\{ b_0,b_1,b_2,b_3 ; c_1,c_2,c_3,c_4 \} = \{ 56,45,16,1;$ $1,8,45,56 \}$ in the following way. Let $C$ be the smallest conjugacy class of elements of order 3 in the group $Suz{:}2$. Let $c$ denote an arbitrary element of $C$ and denote by $V$ the set of all $c' \in C$ for which $cc'$ has order 2. Then $\Upsilon$ is the graph with vertex set $V$, where two vertices $x$ and $y$ are adjacent whenever $xy$ has order 2.

If $x$ is a vertex of $\Upsilon$, then by \cite{So} we know that $|\Upsilon_2(x)|=315$ and that the local graph $\Upsilon_x$ is isomorphic to the Gewirtz graph. The stabilizer $G_x$ of $x$ inside $G = \mathrm{Aut}(\Upsilon)$ is isomorphic to $\mathrm{L}_3(4){:}2^2$ whose natural action on $\Upsilon_1(x)$ is equivalent with the natural action of $\mathrm{Aut}(\Gamma)$ on $\Gamma$.

We have computationally verified (see \cite{ab-bdb:3}) that if $y \in \Upsilon_2(x)$, then the set $\Upsilon_1(x) \cap \Upsilon_1(y)$ is a special 8-set $X_y$ of the Gewirtz graph $\Upsilon_x$, i.e. there exists a unique central involution $\sigma_y$ of $G_x \cong \mathrm{L}_3(4){:}2^2$ fixing each element of $X_y$. The map $y \mapsto \sigma_y$ defines a bijection between $\Upsilon_2(x)$ and the set of 315 central involutions of $G_x \cong \mathrm{L}_3(4){:}2^2$. In view of this bijection, the natural action of $G_x \cong \mathrm{L}_3(4){:}2^2$ on $\Upsilon_2(x)$ is equivalent with the action of $\mathrm{L}_3(4){:}2^2$ on the set of central involutions of $\mathrm{L}_3(4){:}2^2$ (given by conjugation). 
Using the same notation for the suborbits as given in Section \ref{sec2}, we have computationally verified the following for two vertices $y_1,y_2 \in \Upsilon_2(x)$ and $Y = \Upsilon_1(y_1) \cap \Upsilon_1(y_2)$.

\begin{center}
\begin{tabular}{|c|c|c|c|c|}
\hline
 & $\d(y_1,y_2)$ & $|\Upsilon_1(x) \cap Y|$ & $|\Upsilon_2(x) \cap Y|$ & $|\Upsilon_3(x) \cap Y|$ \\
\hline
\hline
$y_2 \in \mathcal{O}_0(y_1)$ & 0 & 8 & 32 & 16 \\
\hline
$y_2 \in \mathcal{O}_{1a}(y_1)$ & 4 & 0 & 0 & 0 \\
\hline
$y_2 \in \mathcal{O}_{1b}(y_1)$ & 2 & 4 & 4 & 0 \\
\hline
$y_2 \in \mathcal{O}_{2a}(y_1)$ & 2 & 0 & 4 & 4 \\
\hline
$y_2 \in \mathcal{O}_{2b}(y_1)$ & 1 & 2 & 4 & 4 \\
\hline
$y_2 \in \mathcal{O}_{3a}(y_1)$ & 3 & 0 & 0 & 0 \\
\hline
$y_2 \in \mathcal{O}_{3b}(y_1)$ & 2 & 2 & 4 & 2 \\
\hline
$y_2 \in \mathcal{O}_4(y_1)$ & 2 & 1 & 4 & 3 \\
\hline
\end{tabular}
\end{center}

\medskip \noindent Theorem \ref{theo1.4} readily follows from the information mentioned in the table.

\appendix
\section{The $\mathrm{G}_2(4)$ near octagon and a distance-regular graph of Soicher}

In this appendix we establish the connection between the near octagon $\O_1$ and the distance-regular graph $\Sigma$ discovered by Soicher in \cite{So}, which has intersection array $\{416, 315, 64, 1; 1, 32, 315, 416\}$. 
The graph $\Sigma$ is a triple cover of the Suzuki graph (in the sense of \cite{So}) and has $\mathrm{Aut}(\Sigma) \cong 3 \cdot Suz{:}2$. 
Under the action of the stabilizer (which is isomorphic to $\mathrm{G}_2(4){:}2$) of its automorphism group with respect to a vertex $x$, the orbits are equal to $\Sigma_0(x)$, $\Sigma_1(x)$, $\Sigma_2(x)$, $\Sigma_3(x)$ and $\Sigma_4(x)$, with sizes $1$, $416$, $4095$, $832$ and $2$, respectively (so in particular, 
the graph $\Sigma$ is distance-transitive).

Let $x$ be a fixed vertex of $\Sigma$.  The local graph $\Sigma_x$ is isomorphic to the well-known $\mathrm{G}_2(4)$-graph, which is a strongly regular graph with parameters $(v, k, \lambda, \mu) = (416, 100, 36, 20)$. Let $G = \mathrm{Aut}(\Sigma)$. Then the stabilizer $G_x$ is isomorphic to the group $\mathrm{G}_2(4){:}2$, and it is the full automorphism group of the local graph $\Sigma_x$. Let $\O_1(x)$ denote the near octagon defined on the central involutions of $G_x \cong \mathrm{G}_2(4){:}2$ (in the sense of Proposition \ref{prop1.1}). We have computationally verified (see \cite{ab-bdb:3}) that every such involution $\sigma$ (which is a point of $\O_1(x)$) fixes $32$ points of the $\mathrm{G}_2(4)$-graph $\Sigma_x$, which we denote by $X_\sigma$. Moreover, the elements of $X_\sigma$ determine the involution $\sigma$ uniquely, as the group generated by $\sigma$ is the unique subgroup of $\mathrm{G}_2(4){:}2$ which fixes $X_\sigma$ pointwise. Therefore, the $4095$ points of $\O_1(x)$ are in bijective correspondence with $4095$ such \textit{special $32$-sets} in the $\mathrm{G}_2(4)$-graph $\Sigma_x$. 

In Soicher's graph $\Sigma$ we can computationally check that for every vertex $y \in \Sigma_2(x)$, the set $\Sigma_1(x) \cap \Sigma_1(y)$ is a \textit{special $32$-set} of the $\mathrm{G}_2(4)$-graph $\Sigma_x$, i.e., there exists a unique central involution $\sigma_y$ of $\mathrm{G}_2(4){:}2$ which fixes this set $X_y$ pointwise. In this manner, we get a map $y \mapsto \sigma_y$ between the $4095$ elements of $\Sigma_2(x)$ and the $4095$ central involutions of $\mathrm{G}_2(4){:}2$ which is bijective. Moreover, computations in the graph $\Sigma$ give us the following information for two points $y_1, y_2 \in \Sigma_2(x)$ with $Y = \Sigma_1(y_1) \cap \Sigma_1(y_2)$, recorded in Table \ref{tab:G24_fixed2}. 

\begin{table}[!htbp]
\begin{center}
\begin{tabular}{|c|c|c|c|c|}
\hline
 & $\dist(y_1,y_2)$ & $|\Sigma_1(x) \cap Y|$ & $|\Sigma_2(x) \cap Y|$ & $|\Sigma_3(x) \cap Y|$ \\
\hline
\hline
$y_2 \in \mathcal{O}_0(y_1)$ & 0 & 32 & 320 & 64 \\
\hline
$y_2 \in \mathcal{O}_{1a}(y_1)$ & 4 & 0 & 0 & 0 \\
\hline
$y_2 \in \mathcal{O}_{1b}(y_1)$ & 2 & 16 & 16 & 0 \\
\hline
$y_2 \in \mathcal{O}_{2a}(y_1)$ & 2 & 0 & 16 & 16 \\
\hline
$y_2 \in \mathcal{O}_{2b}(y_1)$ & 1 & 8 & 76 & 16 \\
\hline
$y_2 \in \mathcal{O}_{3a}(y_1)$ & 3 & 0 & 0 & 0 \\
\hline
$y_2 \in \mathcal{O}_{3b}(y_1)$ & 2 & 5 & 25 & 2 \\
\hline
$y_2 \in \mathcal{O}_4(y_1)$ & 2 & 1 & 25 & 6 \\
\hline
\end{tabular}
\end{center}
\caption{Intersection patterns in Soicher's first graph}
\label{tab:G24_fixed2}
\end{table}

From Table \ref{tab:G24_fixed2}, it follows that for two vertices $y_1, y_2$ in $\Sigma_2(x)$, the involutions $\sigma_{y_1}$ and $\sigma_{y_2}$ are collinear in the near octagon $\O_1(x)$ (which is equivalent to $\sigma_2 \in \m O_{1a}(\sigma_1) \cup \m O_{1b}(\sigma_1)$) if and only if $\dist(y_1, y_2) = 4$ or ($\dist(y_1, y_2) = 2$ and $|\Sigma_1(x) \cap \Sigma_1(y_1) \cap \Sigma_1(y_2)| = 16$). 
Thus we have the following.

\begin{theo} \label{theoA.1}
Let $x$ be a vertex of $\Sigma$. Then the following hold:
\begin{enumerate}
\item[$(1)$] For every vertex $y \in \Sigma_2(x)$, the elementwise stabilizer of $\Sigma_1(x) \cap \Sigma_1(y)$ inside $G_x$ has order $2$ and is generated by a central involution $\sigma_y$ of $G_x \cong \mathrm{G}_2(4){:}2$. Moreover, the map $\theta: y \mapsto \sigma_y$ defines a bijection between $\Sigma_2(x)$ and the set of central involutions of $G_x$.

\item[$(2)$] Let $\Gamma$ denote the graph defined on the set $\Sigma_2(x)$ of vertices at distance $2$ from $x$ in $\Sigma$, by making two vertices $y_1, y_2$ adjacent if and only if $\dist(y_1, y_2) = 4$ or ($\dist(y_1, y_2) = 2$ and $|\Sigma_1(x) \cap \Sigma_1(y_1) \cap \Sigma_1(y_2)| = 16$). Then the map $\theta$ is an isomorphism between $\Gamma$ and the collinearity graph of $\O_1(x)$.

\item[$(3)$] The map $\theta$ also defines an isomorphism between the subgraph of $\Sigma$ induced on $\Sigma_2(x)$ (the second subconstituent) and the graph obtained from the collinearity graph of $\O_1(x)$ by making two vertices (points of the near octagon) adjacent when they are at distance $2$ from each other and have a unique common neighbour. 
\end{enumerate}
\end{theo}

\end{document}